\RequirePackage{amsmath,amssymb}
\documentclass{jfp}

\usepackage{natbib}
\usepackage{code}
\usepackage[svgnames]{xcolor}
\usepackage{stmaryrd}
\usepackage{stfloats}
\usepackage{url}
\usepackage{jms-proof}
\usepackage{mathpartir}
\usepackage{tikz}
\usepackage{tikz-cd}
\usepackage[framemethod=TikZ]{mdframed}
\usepackage{microtype}

\newtheorem{theorem}{Theorem}[section]
\newtheorem{example}[theorem]{Example}
\newtheorem{lemma}[theorem]{Lemma}
\newtheorem{corollary}[theorem]{Corollary}
\newtheorem{proposition}[theorem]{Proposition}
\newtheorem{definition}[theorem]{Definition}
\newtheorem{construction}[theorem]{Construction}
\newtheorem{notation}[theorem]{Notation}
\newtheorem{remark}[theorem]{Remark}
\newtheorem{convention}[theorem]{Convention}

\usepackage{inconsolata}

\usetikzlibrary{shadings,intersections,arrows,patterns}
\usepackage{pgfplots}
\pgfplotsset{compat=1.11}
\usepgfplotslibrary{fillbetween}

\usepackage{macros}

\begin{document}

\title{Higher-Order Functions and Brouwer's Thesis}

\journaltitle{JFP}
\cpr{Cambridge University Press}
\doival{10.1017/xxxxx}

\lefttitle{Higher-Order Functions and Brouwer's Thesis}
\righttitle{Higher-Order Functions and Brouwer's Thesis}

\jnlDoiYr{2021}

\begin{authgrp}
\author{Jonathan Sterling}
\affiliation{Department of Computer Science,\\
        Carnegie Mellon University\\
        (\email{jmsterli@cs.cmu.edu})}
\end{authgrp}

\begin{abstract}
  Extending Mart\'{\i}n Escard\'o's \emph{effectful forcing} technique, we give
  a new proof of a well-known result: Brouwer's monotone bar theorem holds for
  any bar that can be realized by a functional of type $(\Nat\to\Nat)\to\Nat$
  in G\"odel's \SystemT.  Effectful forcing is an elementary alternative to
  standard sheaf-theoretic forcing arguments, using ideas from programming
  languages, including computational effects, monads, the algebra
  interpretation of call-by-name $\lambda$-calculus, and logical relations.

  Our argument proceeds by interpreting \SystemT{} programs as well-founded
  dialogue trees whose nodes branch on a query to an oracle of type
  $\Nat\to\Nat$, lifted to higher type along a call-by-name translation. To
  connect this interpretation to the bar theorem, we then show that Brouwer's
  famous ``mental constructions'' of barhood constitute an invariant form of
  these dialogue trees in which queries to the oracle are made maximally and in
  order.
\end{abstract}

\maketitle

\begin{quote}
  ``The force of \textbf{Church's Law} is that it postulates that all future notions of
  computation will be equivalent in expressive power (measured by definability
  of functions on the natural numbers) to the $\lambda$-calculus. Church's Law
  is therefore a \emph{scientific law} in the same sense as, say, Newton's Law
  of Universal Gravitation, which makes a prediction about all future
  measurements of the acceleration in a gravitational field.'' \citep{harper:2016}
\end{quote}
\bigskip

\section{Introduction}

It is right to rebel! A major theme of Bob Harper's thinking has been the
critique of universal computation as a scientific explanation of functional
programming, which is fundamentally about \emph{higher-order functions} rather
than functions between the natural numbers. Indeed, many sensible models of
computation do not coincide with Turing computability at higher type, a
critical observation for the era of the \emph{socialization of computation}, in which
programs increasingly act as interacting nodes in complex systems.

An early exponent of this perspective was the Dutch topologist L.E.J.\ Brouwer,
who based his \emph{intuitionistic} program of mathematics on a profound
reconstruction of the notion of higher-order function in terms of infinitary
dialogues, rejecting both the classical emphasis on static relations, and the
(proto-)constructivist fixation on finitude and formal definability.

\subsection{Point-set and constructive topology}\label{sec:pointless}

Classical topology is based on point-sets: a topological space is a set of
points $X$ together with a \emph{frame} of open subsets $\Opens{X}\subseteq
\Pow{X}$, i.e.\ a collection of subsets of $X$ closed under unions and
finite intersections. For instance, the real line $\mathbb{R}$ has the set of
real numbers for points, and the opens are the unions of open intervals
$\cup_\alpha (a_\alpha,b_\alpha)$. A continuous function between topological
spaces $\Mor[f]{X}{Y}$ is a function of the underlying point-sets whose inverse
image $\Mor[f^*]{\Pow{Y}}{\Pow{X}}$ sends opens to opens.

Unfortunately, the theory of topological spaces based on point-sets quickly
becomes non-constructive --- meaning that most important theorems of topology
rely crucially on Zermelo's axiom of choice, or at least the law of the
excluded middle.  Mathematical theorems proved \emph{without} these axioms are
of interest because they correspond to \emph{algorithms} that can, in
principle, be executed by a computer: to put it bluntly, a proof that there
exists a number satisfying some property carries within it an algorithm to
compute such a number.

Brouwer was a second-generation pioneer of this idea that mathematical proof is
always a \emph{programming act}, so he rejected both the general law of the
excluded middle as well as the axiom of choice; as a professional topologist,
Brouwer therefore needed to develop a theory of space that did not place
point-sets at its center. The purpose of this section is to explain how such a
concept of space can exist, and to motivate Brouwer's controversial bar thesis.

In the interest of exposition, we take some
historical liberties in Section~\ref{sec:locales} with our explanation of how
to think of \emph{point-free space}, using more modern language than was
available to Brouwer in his lifetime, before turning to Brouwer's account in
Section~\ref{sec:spreads}.
There are by now many references that introduce modern perspectives on
point-free
topology~\citep{johnstone:1982,johnstone:1983,maclane-moerdijk:1992,sambin:2012,vickers:2007,vickers:1989,anel-joyal:2019}.

\subsubsection{The algebraic structure of opens and constructive spaces}\label{sec:locales}

Many aspects of a topological space can be understood in terms of its open
sets; for instance, if $X$ is a Hausdorff space, the point-set of $X$ can be
reconstructed from the set of \emph{completely prime filters}
on $\Opens{X}$ --- the definition of a prime filter here is not as important as
the fact that the point-set can be reconstructed from purely algebraic data on
open sets.

This suggests that one might define a \emph{new} kind of ``point-free space''
generalizing Hausdorff spaces by taking open sets and inverse image functions
as primitive, and deriving point-sets and continuous functions between point-sets from
these.  It turns out that most aspects of point-set topology that require
classical axioms will become fully constructive when re-stated for these new
point-free spaces.

The algebraic structure of the collection of open sets of a topological space
is concentrated in the concept of a \emph{frame}.

\begin{definition}[Frame]
  A \emph{frame} is a partially ordered set $F$ closed under finite meets $U\land V$ and set-indexed
  joins $\bigvee_{\alpha\in S} U_\alpha$, such that meets distribute across
  joins. A morphism of frames, written $\Mor[f^*]{F}{G}$, is a morphism of the
  underlying poset that preserves finite meets and set-indexed joins.
\end{definition}

The morphisms of frames abstract the behavior of the \emph{inverse images}
$\Mor[f^*]{\Opens{Y}}{\Opens{X}}$ of continuous maps of topological spaces
$\Mor[f]{X}{Y}$. We obtain, therefore, a possible notion of \emph{constructive
space} by considering the formal dual of frames --- by ``dual'' we mean that a
constructive space is just a frame, but a morphism between constructive spaces
goes in the opposite direction.\footnote{We use \emph{constructive space} to
refer to what is called a \emph{locale} in the literature.}

\begin{definition}[Constructive space]
  A \emph{constructive space} $X$ is defined by a frame, conventionally called
  $\Opens{X}$; a \emph{continuous map} of constructive spaces $\Mor[f]{X}{Y}$ is given by a
  morphism of frames $\Mor[f^*]{\Opens{Y}}{\Opens{X}}$, called the
  \emph{inverse image}.
\end{definition}

Then, geometrical figures like points are defined synthetically in terms of
continuous maps. The constructive space corresponding to a single point $*$ is given by the
frame $\Opens{*} = \braces{\bot < \top}$; then, one may define a point of an
arbitrary constructive space $X$ to just be a continuous map $\Mor[x]{*}{X}$.
Unraveling definitions, this turns out to be the same as a completely prime
filter in $\Opens{X}$, but it is not important for our purposes to know this
definition.

\subsubsection{Brouwer's presentation of constructive spaces as spreads}\label{sec:spreads}

Intuitionists (like Brouwer and Weyl) and their forbears (like Poincar\'e) were
not only concerned with the law of the excluded middle and the axiom of choice;
they also did not accept the use of impredicative definitions, where the object being defined may implicitly refer to itself. In
particular, the existence of the \emph{powerset} $\Pow{X}$, the set of
\emph{all} subsets of $X$, was considered doubtful. Unfortunately, many aspects
of the theory of constructive spaces are impredicative in just this sense, but it is
possible to present constructive spaces by generators and relations in a predicative way
--- just like the infinite set of natural numbers can be presented by finitely
many generators in a programming language:

\begin{code}
    datatype nat = Z | S of nat
\end{code}

For his purposes in developing intuitionistic analysis, i.e.\ the study of the
continuum from the intuitionistic perspective, Brouwer developed such a
presentation of constructive spaces called \emph{spreads}, in which the
algebraic structure of opens is made concrete in terms of finite paths in an
infinitely broad and infinitely deep rooted tree.

\begin{notation}[Lists]
  If $S$ is a set, then $\List{S}$ is the set of lists of elements of $S$; we
  will use the following notations:
  \begin{enumerate}
    \item $\Nil$ is the empty list.
    \item $\Singl{u}$ is the singleton list given $u:S$.
    \item $\Prepend{\vec{u}}{\vec{v}}$ is the list obtained by appending $\vec{v}$ to the end of $\vec{u}$.
    \item We write $\Len{\vec{u}}$ for the length of the list $\vec{u}$.
    \item We write $\Cons{v}{\vec{u}}$ for $\Prepend{\Singl{v}}{\vec{u}}$, and $\Snoc{\vec{u}}{v}$ for $\Prepend{\vec{u}}{\Singl{v}}$.
    \item We write $\vec{u}\preceq\vec{v}$ to assert that $\vec{u}$ is a prefix of $\vec{v}$.
  \end{enumerate}
\end{notation}

\begin{notation}[Infinite sequences]
  If $S$ is a set, then $\Nat\to S$ is the set of infinite sequences of elements of $S$, which we manipulate using the following notations:
  \begin{enumerate}
    \item We write $\Take{k}{\alpha}$ to mean the list comprised of the first $k$ elements of $\alpha$.
    \item We write $\InOpen{\alpha}{\vec{u}}$ or $\vec{u}\prec\alpha$ to assert that $\vec{u}:\List{S}$ is a finite prefix of $\alpha:\Nat\to S$.
  \end{enumerate}
\end{notation}

\begin{definition}[Spread law]
  A \emph{spread law} is a decidable set $\Law{}$ of lists
  $\vec{u}:\List{\Nat}$ called ``basic observations'',\footnote{The
  decidability of $\vec{u}\in\Law{}$ of course is only meaningful in a setting
  where we have not assumed the law of the excluded middle.} satisfying the
  following axioms that make $\Law{}$ encode an infinitely wide and deep rooted
  tree:
  \begin{enumerate}
    \item The empty sequence $\Nil$ is in $\Law{}$; this ensures that the resulting tree is rooted.
    \item If $\vec{v}\in\Law{}$, then for any prefix $\vec{u}\preceq\vec{v}$, we
      have $\vec{u}\in\Law{}$.
    \item If $\vec{u}\in\Law{}$, then there exists some $k:\Nat$ such that the
      extension $\Snoc{\vec{u}}{k}$ is in $\Law{}$; this ensures that the tree is
      infinitely deep.
  \end{enumerate}
\end{definition}

To view the spread law as an encoding of a tree, think of each
$\vec{u}\in\Law{}$ as a sequence of directions to move from the root of the
tree to one of its nodes. A \emph{spread}, defined below, is
the geometrical object corresponding to the algebraic data of a spread law.

\begin{definition}[Spread]
  A \emph{spread} $\mathbf{S}$ is defined by the data of a spread law $\Law$; one thinks of
  $\mathbf{S}$ as the ``space of infinite paths'' of the tree
  encoded by $\Law$.
\end{definition}

It is useful to switch freely between the intuition of $\vec{u}$ as a node in a
tree and of $\vec{u}$ as an observation about an infinite path in a tree,
substantiated by the fact that each list $\vec{u}$ determines the set of
infinite sequences that it prefixes.
While the data of a spread is the same as the data of a spread law, we define
them separately because a function between spread laws, suitably defined,
should encode the \emph{inverse image} part of a continuous function between
spreads.

\begin{definition}[Points]\label{def:spread-point}
  A point of a spread $\mathbf{S}$ is an infinite path in the tree encoded by
  $\Law$, i.e.\ an infinite sequence $\alpha:\Nat\to\Nat$ such that each finite
  prefix $\vec{u}\prec \alpha$ is in $\Law$, i.e.\ $\Take{k}{\alpha} \in
  \Law$ for all $k:\Nat$.
\end{definition}

A basic observation in a spread is then conceptually similar to a basic open in
a basis for a topology, e.g.\ an open interval on the real line.  From this
perspective, the empty list $\Nil$ represents the entire spread $\mathbf{S}$,
and a general prefix $\vec{u}$ represents the subspace of $\mathbf{S}$ spanned
by points that begin with $\vec{u}$.

\begin{example}[Real numbers]

  Of central importance is the \emph{continuum} $\mathbf{R}$, the spread of
  real numbers. One construction of this spread involves an encoding of
  rational intervals as natural numbers; then a list of natural
  numbers $\vec{u}$ is admitted in $\Law{\mathbf{R}}$ just when its elements
  encode increasingly smaller nested rational intervals each of the
  form $\big[\frac{a}{2^{n-1}}, \frac{a+1}{2^{n-1}}\big]$. A point $\alpha$ in
  the resulting spread $\mathbf{R}$ is an encoding of a real number, and one
  may identify the collection $\mathbb{R}$ of real numbers with
  the equivalence classes of $\mathbf{R}$-points under the suitable
  relation~\citep{van-atten-van-dalen-tieszen:2002}.
\end{example}

\begin{definition}[Fans, or finitary spreads]\label{def:fan}
  A spread whose nodes are all finitarily branching is called a \emph{fan}; an
  example of such a finitary spread is the fan $\mathbf{I}$ representing the
  closed interval $[0,1]\subset \mathbb{R}$. As above, we may write $\mathbb{I}$
  for the appropriate collection of equivalence classes of $\mathbf{I}$-points.
\end{definition}

\subsection{Brouwer's bar thesis: higher-order functions as trees}

In classical topology, an \emph{open cover} over an open $U\in\Opens{X}$ is a
set of opens $\{V_i\in\Opens{X} \mid i:I\}$ such that every point $x\in U$ is contained in one of the
$V_i$. Brouwer defined an analogous notion of cover for spreads,
which he called a \emph{bar}, also in terms of points.

\begin{definition}[Bars in terms of points]
  Let $\mathbf{S}$ be a spread and let $\vec{u}\in\Law$ be a basic observation;
  then a \emph{bar} over $\vec{u}$ is given by a subset
  $Q\subseteq \Law$ satisfying the following axioms:\footnote{We use the
  term \emph{bar} to mean what much of the literature refers to as a \emph{monotone
  bar}. We do not require that membership in the subset $Q\subseteq\Law$ be decidable.}
  \begin{enumerate}

    \item \emph{Monotonicity.} If $\vec{v}$ is in $Q$ and
      $\vec{v}\preceq\vec{w}$, then $\vec{w}$ is in $Q$.

    \item \emph{Covering.} For every point $\InOpen{\alpha}{\vec{u}}$, there exists
      some basic observation $\vec{v}_\alpha\in Q$ such that $\InOpen{\alpha}{\vec{v}_\alpha}$.
  \end{enumerate}

  \textbf{We write $\Covers{\vec{u}}{Q}$ to assert that $Q$ is a bar over
  $\vec{u}$}; colloquially, the \emph{covering} axiom states that every
  infinite path $\alpha$ eventually ``hits the bar'' $Q$, as depicted in
  Fig.~\ref{fig:diagram-bar}. It is useful to observe that the covering axiom is
  equivalent to the assertion that for every $\InOpen{\alpha}{\vec{u}}$, there
  exists $k:\Nat$ such that $\Take{k}{\alpha}\in Q$.
\end{definition}

Given that the empty observation $\Nil$ represents the entire spread
$\mathbf{S}$ (because every point $\alpha$ is prefixed by $\Nil$), one refers
to a bar of $\mathbf{S}$ as a subset $Q\subseteq\Law$ such that
$\Covers{\Nil}{Q}$.

\begin{figure}
  \centering

\begin{tikzpicture}
  \tikzstyle{every node}=[font=\LARGE]
  \pgfdeclarelayer{pre main}
  \pgfsetlayers{pre main,main}

  \coordinate (O) at (0,0);

  \draw[gray, densely dashed, name path=Q] (-3,4.0) .. controls (-1,5) and (1,1.5) .. (3,4);
  \path[name path=X] (-3,6) to (3,6);

  \draw[->, red, thick, name path=alpha] plot [smooth,tension=1] coordinates {(0,0) (.2,2) (-.5,4) (0,6)};

  \draw[name path=S] plot coordinates {(-2.66,4.66) (0,0) (2.66,4.66)};
  \tikzfillbetween[of=X and Q]{color=Gray!10}

  \coordinate [name intersections={of=Q and alpha, by={K}}];
  \coordinate [label={left:$\mathbf{S}$}] (Slbl) at (1.6,1);
  \coordinate [label={left:$\alpha$}] (alphalbl) at (.8, 2.4);
  \coordinate [label={left:$Q$}] (Qlbl) at (1.4, 4.8);

  \node[circle,fill,inner sep=2pt,outer sep=0pt,label={[font=\large]right:$\vec{v}_\alpha$}] at (K) {};
\end{tikzpicture}

   \caption{A visualization of when a monotone subset
    $\IsSubsetEq{Q}{\Law}$ is a bar, i.e.\
    $\Covers{\Nil}{Q}$. The dotted line indicates the beginning of the bar, and $\alpha$ is an arbitrary point in the
    spread; $\vec{v}_\alpha$ is a prefix of $\alpha$
    lying in $Q$.}\label{fig:diagram-bar}
\end{figure}

\subsubsection{Intuitionistically, bars are higher-order functions}

From the point of view of the Brouwer--Heyting--Kolmogorov interpretation of the
logical connectives~\citep{heyting:1956}, to say that $Q$ is a bar over
$\vec{u}$ is the same as to require a \emph{higher-order function}
$\phi : \parens{\Nat\to\Nat}\to\Nat$ such that for all $\InOpen{\alpha}{\vec{u}}$, we
have $\Take{\phi(\alpha)}{\alpha}\in Q$. Intuitively, this higher order
function receives a point and responds with how far one must proceed along the point's approximations in order
to reach the bar.

\subsubsection{Bars represented by well-founded trees}

It is an empirical reality that one cannot define a discontinuous function on
the reals by computational means --- this is because such functions always rely
on the ability to branch on whether two real numbers are equal, or something
analogous, which cannot be achieved in finite time. Consequently, it was
important for Brouwer's intuitionistic version of analysis to reflect this limitation
in a \emph{continuity theorem}: all functions defined on the reals are
continuous.

In fact, the crown jewel of Brouwer's intuitionistic analysis is the
verification of the even stronger \emph{uniform continuity} theorem for all
total functions defined on the closed interval, which exhibits not only a
modulus of continuity for such functions, but one that is independent of the
input:
\begin{equation}
  \forall F:\mathbb{I}\to\mathbb{R}.\
  \forall \epsilon:\mathbb{R}^+.\
  \exists \delta:\mathbb{R}^+.\
  \forall x_0,x_1 : \mathbb{I}.\
  \vert x_0 - x_1\vert < \delta
  \Rightarrow
  \vert f(x_0) - f(x_1)\vert < \epsilon
  \tag{UC}
\end{equation}

For intuitionistic mathematicians, uniform continuity is important because it
reflects the limitations of human facilities for manipulating infinitary data:
a human can only make finitely many observations of infinite data. The
uniformity condition expresses that, in the case of the closed interval, one
may bound the number of observations \emph{regardless} of the point being
observed. From a computer science point of view, the uniform continuity theorem
suggests a realistic invariant on nodes that consume and emit bits of data:
there is a bound on how many messages it takes to engender a response.

Because Brouwer represented the closed interval $\mathbb{I}$ in terms of the
points of a fan (Definition~\ref{def:fan}), he could reduce the statement of
uniform continuity to a combination of \emph{ordinary} continuity and his
famous fan theorem, a statement about \emph{bars} $Q$ on the fan
$\mathbf{F}$ whose basic observations are all lists of bits $0,1$:
\begin{equation}
  \Covers{\Nil}{Q} \Rightarrow \exists k:\Nat.\ \forall \alpha:\Nat\to\Bit.\ \Take{k}{\alpha}\in Q
  \tag{FT}
\end{equation}

Observe that the fan theorem is nothing more than a quantifier rotation: the
assumption that $Q$ is a bar associates to each sequence $\alpha$ the length
$k$ of a basic observation of $\alpha$ that lands in the bar; the conclusion of
the fan theorem establishes a uniform bound on the length of this finite prefix
for all $\alpha$.

Brouwer argued for his fan theorem by reflecting on the possible ways to
construct a proof of the antecedent $\Covers{\vec{u}}{Q}$; he hypothesized that
such a proof could always be presented as a well-founded tree governed by the following primitive inferences:
\begin{enumerate}
  \item \emph{\boldeta-inference.} If $\vec{u}\in Q$, then $Q$ bars $\vec{u}$.
  \item \emph{\bolddigamma-inference.}\footnote{Brouwer used \bolddigamma, an archaic Greek letter pronounced ``digamma''.} If $Q$ bars every immediate subnode $\Snoc{\vec{u}}{x}$ of $\vec{u}$, then $Q$ bars $\vec{u}$.
\end{enumerate}

In other words, Brouwer believed that a proof of $\Covers{\vec{u}}{Q}$ could
always be tabulated into a proof of an alternative, \emph{inductive}
characterization of barhood $\IndCovers{\vec{u}}{Q}$, defined relative to a spread $\mathbf{S}$:
\begin{mathpar}
  \infer[\boldeta]{
    \IndCovers{\vec{u}}{Q}
  }{
    \vec{u}\in Q
  }
  \and
  \infer[\bolddigamma]{
    \IndCovers{\vec{u}}{Q}
  }{
    \vec{u}\in\Law
    &
    \forall x:\Nat.\
    \Snoc{\vec{u}}{x}\in\Law{\mathbf{S}}\Rightarrow
    \IndCovers{\Snoc{\vec{u}}{x}}{Q}
  }
\end{mathpar}

\begin{figure}
  \centering

\begin{tikzpicture}
  [level distance=10mm,
   every node/.style={fill=MidnightBlue!20, draw=Black!30, circle, inner sep=1pt, font=\large},
   level 1/.style={sibling distance=45mm,nodes={fill=MidnightBlue!15}},
   level 2/.style={sibling distance=15mm,nodes={fill=MidnightBlue!10}},
   level 3/.style={sibling distance=15mm,nodes={fill=MidnightBlue!5}},
   level 4/.style={sibling distance=15mm,nodes={fill=MidnightBlue!3}}]

  \pgfdeclarelayer{pre main}
  \pgfsetlayers{pre main,main}

  \node  (Nil) {$\bolddigamma_{\Nil}$} [grow'=up]
  child {node (O) {$\bolddigamma_{\mathtt{0}}$}
    child {node (OO) {$\boldeta_{\mathtt{00}}$}
    }
    child {node {$\bolddigamma_{\mathtt{01}}$}
      child {node (OIO) {$\boldeta_{\mathtt{010}}$}}
      child {node (OII) {$\boldeta_{\mathtt{011}}$}}
    }
  }
  child {node (I) {$\bolddigamma_{\mathtt{1}}$}
    child {node (IO) {$\bolddigamma_{\mathtt{10}}$}
      child {node (IOO) {$\bolddigamma_{\mathtt{100}}$}
        child {node (IOOO) {$\boldeta_{\mathtt{1000}}$}}
        child {node (IOOI) {$\boldeta_{\mathtt{1001}}$}}
      }
      child {node (IOI) {$\boldeta_{\mathtt{101}}$}}
    }
    child {node (II) {$\boldeta_{\mathtt{11}}$}}
  };

  \path[name path=X] (-4.5,6.5) to (4.5,6.5);
  \draw[black, thick, name path=L] (Nil) -- (O) -- (OO) -- (-4.5,6.5);
  \draw[black, thick, name path=R] (Nil) -- (I) -- (II) -- (4.5,6.5);
  \draw[gray, thick, densely dashed, name path=Q, smooth] (OO) -- (OIO) -- (OII) -- (IOOO) -- (IOOI) -- (IOI) -- (II);
  \path[name path=Q'] (OO.center) -- (OIO) -- (OII) -- (IOOO) -- (IOOI) -- (IOI) -- (II.center);

  \draw[red, ->, line width=1.5pt] (Nil) -- (I) -- (IO) -- (IOO) -- (IOOI) -- (-.5,6.6);

  \coordinate [label={[font=\LARGE]left:$Q$}] (Qlbl) at (3.3, 5);
  \coordinate [label={[font=\LARGE]left:$\mathbf{F}$}] (Slbl) at (4.2,3.3);
  \coordinate [label={[font=\LARGE]left:$\alpha$}] (alphalbl) at (0.8, 5.8);

  \tikzfillbetween[of=Q' and X]{color=Gray!10}

\end{tikzpicture}

   \caption{A visualization of the \emph{inductive} characterization of
    barhood $\IndCovers{\Nil}{Q}$ in a fan $\mathbf{F}$ consisting of sequences of bits. (Contrast with
    Fig.~\ref{fig:diagram-bar}.) In the depicted tree, each node is
    labeled with the list of bits that it codes; to be
    precise, a node labeled $\boldeta_{\vec{u}}$ represents a
    proof of $\IndCovers{\vec{u}}{Q}$ that is a leaf node, and
    a node labeled $\bolddigamma_{\vec{u}}$ represents a proof of
    $\IndCovers{\vec{u}}{Q}$ that is a branching
    node.}\label{fig:diagram-inductive-bar}
\end{figure}

One way to think of these trees is as a \emph{dialogue} between an actor and a
choice sequence $\alpha$: at each step, the actor has either decided upon an
answer (the \boldeta-inference), or asks the choice sequence what its next
element is (the \bolddigamma-inference), and continues to interact using that
information.

If every proof of $\Covers{\vec{u}}{Q}$ could in fact be tabulated into such a
dialogue $\IndCovers{\vec{u}}{Q}$ (presented visually in
Fig.~\ref{fig:diagram-inductive-bar}), then it would be a simple matter to
establish the uniform bound $k$ by induction:
\begin{enumerate}
  \item \emph{\boldeta-inference.} We have $\vec{u}\in Q$, so assign $k =\Len{\vec{u}}$.
  \item \emph{\bolddigamma-inference.} By induction, we have for all possible extensions $\Snoc{\vec{u}}{x}$, a bound $k_x$; because $\mathbf{F}$ is a fan, there are only finitely many such extensions $\Snoc{\vec{u}}{x}$, so we simply choose $k = \max\{k_x\mid x\}$.
\end{enumerate}

\begin{theorem}[Soundness of inductive barhood]\label{thm:soundness}
  The inductive characterization of barhood is sound: if $\IndCovers{\vec{u}}{Q}$ is proved, then $\Covers{\vec{u}}{Q}$ holds.
\end{theorem}
\begin{proof}
  We proceed by induction on the proof of $\IndCovers{\vec{u}}{Q}$ to prove that
  prove that $\Covers{\vec{u}}{Q}$.
  \begin{enumerate}
    \item \emph{\boldeta-inference.} Then $\vec{u}\in Q$, so we choose $k =\Len{\vec{u}}$.
    \item \emph{\bolddigamma-inference.}
      By induction, we have $\Covers{\Snoc{\vec{u}}{x}}{Q}$ for all possible extensions $x$ with a bound $k$; but this implies $\Covers{\vec{u}}{Q}$ immediately, with bound $k+1$.\qedhere
  \end{enumerate}
\end{proof}

In order to prove his fan theorem (and thence establish uniform continuity for
real functions on the closed interval), Brouwer wished to prove a converse to
Theorem~\ref{thm:soundness}:

\begin{proposition}[Brouwer's bar thesis]\label{prop:completeness}
  The inductive characterization of barhood is complete: if $\Covers{\vec{u}}{Q}$ then $\IndCovers{\vec{u}}{Q}$.
\end{proposition}

\subsection{Functional programming interpretation: exceptions as shared secrets}

It is possible to give an interpretation of Brouwer's bar thesis in terms
of effectful functional programming in Standard~ML, in the spirit of the
``seemingly impossible functional programs'' of \citet{escardo:blog:seemingly-impossible-functional-programs,bauer:blog:sometimes-all-functions-are-continuous,rahli-bickford-constable:2017}.  As we noted, a bar
$\Covers{\vec{u}}{Q}$ can be thought of as an operation that assigns to each
infinite sequence $\InOpen{\alpha}{\vec{u}}$ a number $k:\Nat$ such that
$\Take{k}{\alpha}\in Q$.  Hence, we may represent such ``functional bars'' by a
type of higher-order functions:

\begin{code}
    type fun\_bar = (nat $\to$ nat) $\to$ nat
\end{code}

On the other hand, an inductive bar can be represented in Standard~ML by the
following datatype of infinitely branching trees, in which the \textcd{SPIT} constructor
encodes the \boldeta-inference and the \textcd{BITE} constructor encodes the
\bolddigamma-inference:

\begin{code}
    datatype ind\_bar =
       SPIT of nat
     | BITE of (nat $\to$ ind\_bar)
\end{code}

The argument of \textcd{SPIT} contains the number of steps that were necessary to reach the bar.
Then, we may define an (effectful) algorithm to tabulate a functional bar \textcd{Q
: fun\_bar} as a tree \textcd{tabulate Q : ind\_bar}, keeping track of the list of
numbers \textcd{xs : nat list} we have received from \textcd{BITE} so far. The
tabulation loop works in the following way:
\begin{enumerate}

  \item If \textcd{Q $\alpha$ = Q $\beta$} gives the same answer for all \textcd{$\alpha$,$\beta$ : nat $\to$ nat} that start with \textcd{xs}, then we have reached the bar; therefore return \textcd{SPIT}.

  \item Otherwise, we require at least another digit of information; call
    \textcd{BITE}, abstracting \textcd{x:nat} and recursing with \textcd{xs := xs @ [x]}.

\end{enumerate}

The first step is non-trivial, and we may accomplish it using a computational
effect. Given a fresh exception \textcd{E : exn}, we may define the ``generic
infinite sequence'' extending \textcd{xs} by indexing into \textcd{xs} and raising
\textcd{E} when out of bounds:

\begin{code}
    fun generic (E : exn) (xs : nat list) : nat $\to$ nat =
      fn i $\Rightarrow$
        List.nth (xs, i) handle
          Subscript $\Rightarrow$ raise E
\end{code}

The idea is that we may use the exception \textcd{E} as a \emph{shared secret} in
our tabulation loop\footnote{The view of Standard~ML--style exceptions as
shared secrets is explored by Harper in his blog and in his
book~\citep{harper:blog:exceptions-are-shared-secrets,harper:2016}.}, to
discover that \textcd{Q} has not yet converged and requires a further digit.
\begin{code}
    fun tabulate (Q : fun\_bar) : ind\_bar =
      let
        fun loop xs =
          let
            exception Secret
          in
            SPIT (Q (generic Secret xs)) handle
              Secret $\Rightarrow$ BITE (fn x $\Rightarrow$ loop (xs @ [x]))
          end
      in
        loop []
      end
\end{code}

The use of a computational effect in the Standard~ML implementation of the
tabulation is the essence of what is meant by Escard\'o's term \emph{effectful
forcing}~\citep{escardo:2013}; indeed, the proof of our main result follows the
moral arc of the code presented above, likewise hinging on the existence of a
generic sequence that has the ability to ``spy'' on the internals of any
function that calls it in the manner of a debug harness.

\subsubsection{Relation to the bar thesis} Does the \textcd{tabulate} function
terminate? The answer to this question depends entirely on what arguments \textcd{Q
: fun\_bar} actually can be defined, which is in essence what Brouwer's bar
thesis aims to control. Constructive mathematics is compatible with
interpretations in which certain bizarre higher-order functions \textcd{Q} exist
for which \textcd{tabulate Q} does not terminate; the most famous such case arises
from the \emph{Kleene tree}, an object that exists in the Church--Turing
interpretation of second-order functions as first-order functions on
codes~\citep{troelstra-vandalen:1988}.

\subsection{Perspective and contribution}

Returning to the more modern perspective on constructive topology presented in
Section~\ref{sec:pointless}, the bar thesis enunciated in
Proposition~\ref{prop:completeness} states that a certain constructive space
$\mathbf{B}$ is completely determined by the topological space obtained from
its points; but the spatiality of $\mathbf{B}$, a classical
triviality~\citep{dummett:elements}, is independent of constructive mathematics.

Brouwer's ``proof'' of Proposition~\ref{prop:completeness} was therefore not
successful.  Considering the Brouwer--Heyting--Kolmogorov interpretation of a
proof of $\Covers{\vec{u}}{Q}$ as a higher-order function, to accept Brouwer's
bar thesis is to require that higher-order functions of a certain type can
always be tabulated into well-founded trees or dialogues.  Under the standard
Church--Turing interpretation of higher-order functions, this tabulation is
refuted~\citep{troelstra-vandalen:1988}; far from a no-go theorem, this is in
fact an invitation to consider alternative notions of computability that
justify Brouwer's identification of higher-order functions
$\parens{\Nat\to\Nat}\to\Nat$ with inductive dialogues.

In this paper, we explain one possible notion of computability under which
Brouwer's bar thesis is justified: definability in G\"odel's \SystemT. Our work
is an extension of Mart\'in Escard\'o's \emph{effectful forcing} technique
\citep{escardo:2013}, including a novel normalization theorem for certain interaction
trees into Brouwer's $(\boldeta,\bolddigamma)$-trees.

\subsection{Formalization in Agda}

The main results in this paper have been formalized in the Agda proof
assistant~\citep{norell:2009}. We carried out our proof in Martin-L\"of's
Intensional Type Theory extended by the function extensionality axiom, which
was needed in order to prove the monad laws for Escard\'o's dialogue monad; we
do not assume the unicity of identity proofs, and therefore our proof could in
principle be converted to Cubical Agda~\citep{vezzosi-mortberg-abel:2019}, in
which function extensionality computes.
Our formalization can be accessed as follows:

\begin{enumerate}
  \item Browsable source: \url{http://jonsterling.github.io/agda-effectful-forcing}
  \item Repository: \url{http://www.github.com/jonsterling/agda-effectful-forcing}
\end{enumerate}

\section{G\"odel's \texorpdfstring{\SystemT{}}{System T}}

In his famous \emph{Dialectica} paper, Kurt G\"odel
introduced \SystemT{} to serve as a formal theory of constructions for Heyting
arithmetic~\citep{goedel:1958}; in modern terms, \SystemT*{} is a minimal, total functional
programming language with function types $\TArr{\sigma}{\tau}$ and a base type
of natural numbers $\TNat$.

\begin{grammar}
  contexts & \Gamma,\Delta & \cdot \GrmSep \Gamma,x:\sigma
  \\
  types & \sigma,\tau & \TNat \GrmSep \TArr{\sigma}{\tau}
  \\
  terms & M,N & x \GrmSep \TLam{x}{M} \GrmSep \TAp{\sigma}{M}{N} \GrmSep
  \TZe \GrmSep \TSu{M} \GrmSep
  \\
  \GrmContinue &
  \TRec{N}{M_z}{x}{h}{M_s}
\end{grammar}

In the term for functional abstraction $\TLam{x}{M}$, the variable $x$ is
bound; likewise, in the term for the recursor $\TRec{N}{M_z}{x}{h}{M_s}$, the
variables $x,h$ are both bound.

\begin{convention}[Variable binding]
  In everything that follows, terms are considered modulo $\alpha$-equivalence
  (i.e.\ up to permutations of bound variables), and we assume the ability to
  freely obtain a fresh variable that is not free in a term.
\end{convention}

\subsection{Static semantics of \texorpdfstring{\SystemT{}}{System T}}

The theory of \SystemT{} is given by two forms of judgment:
\begin{enumerate}
  \item $\boxed{\IsCtx{\Gamma}}$ asserts that $\Gamma$ is a well-formed context.
    Contexts assign types to the free variables that are in scope.

  \item $\boxed{\IsTm{\Gamma}{M}{\sigma}}$ asserts that $M$ is a well-formed element of type $\sigma$ in context $\Gamma$; this judgment is only defined when $\IsCtx{\Gamma}$.
\end{enumerate}

We omit definitional equality entirely, as our dialogue model will be too
intensional to validate the usual equations of \SystemT{}.\footnote{We thank
Pierre-Marie P\'edrot for pointing this out to us.}
We do not need a separate judgment $\sigma\ \mathit{type}$, because all raw types
$\sigma$ generated by the grammar above are trivially well-formed.
The rules for forming contexts are quite simple, and have the main role of
preventing variable shadowing (a simplifying measure that will make semantic
interpretation easier):
\begin{mathpar}
  \inferrule{
  }{
    \IsCtx{\cdot}
  }
  \and
  \inferrule{
    \IsCtx{\Gamma}
    \\
    x\not\in\Dom{\Gamma}
  }{
    \IsCtx{\Gamma,x:\sigma}
  }
  \\
  \mbox{(where $\Dom{\Gamma}$ is the set of variables declared in $\Gamma$)}
\end{mathpar}

The well-typed terms are generated by standard rules of inference; for
instance, a variable $x$ is well-formed if it appears in the context:
\begin{mathpar}
  \inferrule[variable]{
  }{
    \IsTm{\Gamma,x:\sigma,\Delta}{x}{\sigma}
  }
\end{mathpar}

\subsubsection{The function type}

The elements of the function type $\TArr{\sigma}{\tau}$ are characterized by an
abstraction and an application rule. In our version of \SystemT*{}, we have
included just enough typing annotations in the syntax to guarantee the
existence of a coherent interpretation of \emph{judgments}
$\IsTm{\Gamma}{M}{\sigma}$ into models; for this reason, we include a typing
annotation on function application $\TAp{\sigma}{M}{N}$, but have no need for a
typing annotation on abstraction $\TLam{x}{M}$:
\begin{mathpar}
  \inferrule[function abstraction]{
    \IsTm{\Gamma,x:\sigma}{M}{\tau}
  }{
    \IsTm{\Gamma}{\TLam{x}{M}}{\TArr{\sigma}{\tau}}
  }
  \and
  \inferrule[function application]{
    \IsTm{\Gamma}{M}{\TArr{\sigma}{\tau}}
    \\
    \IsTm{\Gamma}{N}{\sigma}
  }{
    \IsTm{\Gamma}{\TAp{\sigma}{M}{N}}{\tau}
  }
\end{mathpar}

\subsubsection{Natural numbers and recursion}

The natural numbers are given a unary encoding in our presentation of
\SystemT*{}, where $\TZe$ codes $0$ and $\TSu{N}$ codes $n+1$ if $N$ codes $n$;
a numeral $n:\Nat$ is therefore coded as the iterated application
$\TKwd{s}^n(\TZe)$.
\begin{mathpar}
  \inferrule[zero]{
  }{
    \IsTm{\Gamma}{\TZe}{\TNat}
  }
  \and
  \inferrule[successor]{
    \IsTm{\Gamma}{M}{\TNat}
  }{
    \IsTm{\Gamma}{\TSu{M}}{\TNat}
  }
\end{mathpar}

Programs on natural numbers in \SystemT*{} are written by \emph{primitive
recursion}, at possibly higher type.\footnote{Sometimes the term
\emph{primitive recursion} is restricted to mean recursion at base type.} The
recursor is typed as follows:
\begin{mathpar}
  \inferrule[recursion]{
    \IsTm{\Gamma}{N}{\TNat}
    \\
    \IsTm{\Gamma}{M_z}{\sigma}
    \\
    \IsTm{\Gamma,x:\TNat,h:\sigma}{M_s}{\sigma}
  }{
    \IsTm{\Gamma}{
      \TRec{N}{M_z}{x}{h}{M_s}
    }{\sigma}
  }
\end{mathpar}
The recursor should be thought of as a construct for defining a function
$F:\TArr{\TNat}{\sigma}$ and then immediately applying it to the target $N$; in
the case for the successor, one binds not only $x:\TNat$ but also $h:\sigma$,
which intuitively stands for $F(x)$ (the induction hypothesis).

\subsection{Programming in \texorpdfstring{\SystemT{}}{System T}}

To gain some intuition for higher-order primitive recursion, we work through some examples.

\begin{convention}[Informal notation for \SystemT*{}]
  In order to make our examples more clear, we impose some informal notation for terms.
  \begin{enumerate}

    \item Nested function abstractions are written all at once: for instance
      $\TLam{x,y}{M}$ means $\TLam{x}{\TLam{y}{M}}$.

    \item Numeric constants $\bar{n}$ can be written instead of their
      encodings; for instance $\bar{2}$ means $\TSu{\TSu{\TZe}}$.

    \item Function applications can be written $M(N)$ rather than
      $\TAp{\sigma}{M}{N}$, when the $\sigma$ can be easily inferred from
      context.

    \item A function defined by primitive recursion can be written as
      $\TRec*{M_z}{x}{h}{M_s}$ rather than the abstracted
      $\TLam{y}{\TRec{y}{M_z}{x}{h}{M_s}}$.

  \end{enumerate}
\end{convention}

The simplest example of programming using primitive recursion is given by the
addition function:
\begin{align*}
  \TKwd{plus} &: \TArr{\TNat}{\TArr{\TNat}{\TNat}}
  \\
  \TKwd{plus} &\triangleq
  \TLam{x,y}{
    \TRec{x}{y}{\_}{h}{\TSu{h}}
  }
\end{align*}

Another version of the addition function can be given that uses primitive
recursion at higher type:
\begin{align*}
  \TKwd{plus}' &: \TArr{\TNat}{\TArr{\TNat}{\TNat}}
  \\
  \TKwd{plus}' &\triangleq
  \TRec*{\parens*{\TLam{y}{y}}}{\_}{h}{
    \parens*{
      \TLam{y}{\TSu{h(y)}}
    }
  }
\end{align*}

Higher-order primitive recursion is stronger than ordinary primitive recursion;
for instance, Ackermann's function $\mathcal{A}:\Nat\times\Nat\to\Nat$ is not
primitive recursive, but we will see that it can be defined in
\SystemT*{}.\footnote{A three-parameter variant of Ackermann's function was
discussed by Hilbert in his 1925 lecture \emph{\"Uber das
Unendliche}~\citep{hilbert:1926}, later proved by Ackermann to not be
primitive recursive~\citep{ackermann:1928}.}

\begin{align*}
  \mathcal{A}(0,n) &= n + 1 \\
  \mathcal{A}(m+1,0) &= \mathcal{A}(m,1)\\
  \mathcal{A}(m+1,n+1) &= \mathcal{A}(m,\mathcal{A}(m+1,n))
\end{align*}

In \SystemT*{}, Ackermann's function $\mathcal{A}$ is encoded by a program
$\TKwd{ack}:\TArr{\TNat}{\TArr{\TNat}{\TNat}}$ that computes a \emph{function}
by recursion, as in Harper~\citep{harper:2016}:
\begin{align*}
  \TKwd{ack} &: \TArr{\TNat}{\TArr{\TNat}{\TNat}}
  \\
  \TKwd{ack} &\triangleq
  \TRec*{\parens*{\TLam{n}{\TSu{n}}}}{\_}{g}{
    \parens*{
      \TLam{n}{
        \TKwd{iter}(g)(n)(g(\bar{1}))
      }
    }
  }
\end{align*}
where we define the $n$-fold iteration of a function as follows:
\begin{align*}
  \TKwd{iter} &: \TArr{\parens*{\TArr{\TNat}{\TNat}}}{\TArr{\TNat}{\parens*{\TArr{\TNat}{\TNat}}}}
  \\
  \TKwd{iter} &=
  \TLam{f}{
    \TRec*{\parens*{\TLam{x}{x}}}{\_}{g}{\parens*{\TLam{x}{f(g(x))}}}
  }
\end{align*}

The ability to use recursion to define functions of arbitrary order is the
hallmark of \SystemT{}.

\subsection{Standard semantics of \texorpdfstring{\SystemT{}}{System T}}

In order to state what it means for a bar $Q$ to be ``realizable in
\SystemT{}'', we must first explain the standard/Tarski-style semantics of
\SystemT{}, in which $\TNat$ is interpreted by the ordinary set of natural
numbers. We begin by defining the basic domains of interpretation, starting
with contexts and types:
\begin{align*}
  \StdDomG{\Gamma} &= \textstyle\prod_{x\in \Dom{\Gamma}} \StdDomV{\Gamma(x)}
  \\
  \StdDomV{\TNat} &= \Nat
  \\
  \StdDomV{\TArr{\sigma}{\tau}} &= \StdDomV{\sigma}\to {\StdDomV{\tau}}
\end{align*}

A term $\IsTm{\Gamma}{M}{\sigma}$ is interpreted as a \emph{function} from
$\StdDomG{\Gamma}$ to $\StdDomV{\sigma}$ in the obvious way:
\begin{align*}
  \StdInterp{\Gamma}{M}{\sigma} &: \StdDomG{\Gamma}\to\StdDomV{\sigma}
  \\
  \StdInterp{\Gamma,x:\sigma,\Delta}{x}{\sigma}g &= g(x)
  \\
  \StdInterp{\Gamma}{\TLam{x}{M}}{\TArr{\sigma}{\tau}}g &=
  \lambda a.\ \StdInterp{\Gamma,x:\sigma}{M}{\tau}(g,x\mapsto a)
  \\
  \StdInterp{\Gamma}{\TAp{\sigma}{M}{N}}{\tau}g &=
  \parens*{\StdInterp{\Gamma}{M}{\TArr{\sigma}{\tau}}g}\parens*{\StdInterp{\Gamma}{N}{\sigma}g}
  \\
  \StdInterp{\Gamma}{\TZe}{\TNat}g &= 0
  \\
  \StdInterp{\Gamma}{\TSu{M}}{\TNat}g &= 1 + \StdInterp{\Gamma}{M}{\TNat}g
  \\
  \StdInterp{\Gamma}{\TRec{N}{M_z}{x}{h}{M_s}}{\sigma}g &=
  \Rec_\sigma
    \begin{pmatrix*}[l]
      \StdInterp{\Gamma}{N}{\TNat}g,\\
      \StdInterp{\Gamma}{M_z}{\sigma}g, \\
      \lambda (a,b).\ \StdInterp{\Gamma,x:\TNat,h:\sigma}{M_s}{\sigma}(g,x\mapsto a,h\mapsto b)
    \end{pmatrix*}
\end{align*}
where we define a meta-level recursor as follows:
\begin{align*}
  \Rec_\sigma &: \Nat\times \StdDomV{\sigma}\times
  \parens{
    {\Nat\times\StdDomV{\sigma}}\to
    \StdDomV{\sigma}
  }
  \to \StdDomV{\sigma}\\
  \Rec_\sigma\parens*{0, z, s} &= z\\
  \Rec_\sigma\parens*{n+1, z, s} &=
  s\parens*{n, \Rec_\sigma\parens*{n,z,s}}
\end{align*}

\subsubsection{\texorpdfstring{\SystemT{}}{System T}-realizable bars}

Every closed term $\IsTm{\cdot}{F}{\TArr*{\TArr{\TNat}{\TNat}}{\TNat}}$ determines, via the
standard semantics above, an assignment of natural numbers
$\TMetaAp{F}{\alpha}:\Nat$ to each sequence $\alpha:\Nat\to\Nat$:
\[
  \TMetaAp{F}{\alpha} \triangleq
  \StdInterp{\cdot}{F}{\TArr*{\TArr{\TNat}{\TNat}}{\TNat}}(\alpha)
\]%

Using this notation, we are equipped to define a variation on $\Covers{\vec{u}}{Q}$ that
expresses realizability in \SystemT{}, which we will write
$\TCovers{\vec{u}}{Q}$:
\[
  \parens*{\TCovers{\vec{u}}{Q}}
  \triangleq
  \exists F.\
  \parens*{\IsTm{\cdot}{F}{\TArr*{\TArr{\TNat}{\TNat}}{\TNat}}}
  \land
  \forall\alpha:\Nat\to\Nat.\
  \Prepend{\vec{u}}{
    \Take{\TMetaAp{F}{\alpha}}{\alpha}
  }
  \in Q
\]

\subsection{Interactive semantics of \texorpdfstring{\SystemT{}}{System T}}

Mart\'{\i}n Escard\'o pioneered a technique called ``effectful
forcing'' for demonstrating non-constructive (Brouwerian) principles for the definable
functionals of G\"odel's \SystemT~\citep{escardo:2013}, including the continuity
of functionals $(\Nat\to\Nat)\to\Nat$ and uniform continuity of functionals
$(\Nat\to2)\to\Nat$.

The idea of effectful forcing is to give a non-standard \emph{interactive} model of \SystemT{} in
which an element of type $\TNat$ is interpreted not as a single natural number,
but as a kind of infinitely-branching well-founded tree with natural numbers at the leaves;
the branch nodes of the tree represent queries to an oracle (a choice
sequence). Escard\'o used this model to show that definable functionals in
\SystemT{} are always continuous; our contribution is to extend Escard\'o's
method to prove that the realizable bars can be tabulated into Brouwer's
$(\boldeta,\bolddigamma)$-trees, establishing a restricted version of the bar theorem.

Intuitively, we would to arrange our interactive model in such a way that the
naturals are interpreted quite literally as Brouwer's
$(\boldeta,\bolddigamma)$-trees, but it will be more clear to \emph{first}
interpret \SystemT{} into a more flexible shape of tree (which we will
call ``Escard\'o dialogues''), in which queries to the choice sequence can be
made out of order and repeated. Then, we prove that Brouwer's
$(\boldeta,\bolddigamma)$-trees can be construed as a \emph{normal} form for
the Escard\'o dialogues.

\subsubsection{Escard\'o dialogues: ideal codes for functionals}

Escard\'o's dialogue trees provide, for any sets $I,J,A\in\Sets$, a
representation $\ETree{I}{J}{A}\in\Sets$ of computations of an element of $A$
relative to an \emph{oracle} of type $I\to J$ (i.e.\ an oracle that receives
queries of type $I$ and responds with answers of type $J$). The Escard\'o
dialogues are generated by the following rules of formation:

\begin{mathpar}
  \boxed{
    \inferrule{
      I,J,A : \Sets
    }{
      \ETree{I}{J}{A} : \Sets
    }
  }
  \and
  \inferrule[\RuleERet]{
    a : A
  }{
    \ERet{a} : \ETree{I}{J}{A}
  }
  \and
  \inferrule[\RuleEQuery]{
    i : I
    \\
    \mathfrak{a}: \parens{J\to \ETree{I}{J}{A}}
  }{
    \EQuery{i}{\mathfrak{a}}\in\ETree{I}{J}{A}
  }
\end{mathpar}

The Escard\'o dialogues exhibit an endofunctor $\Sets\to\Sets$ in the following
way:
\begin{align*}
  \EMap&: \parens{A\to B}\to \ETree{I}{J}{A}\to\ETree{I}{J}{B}
  \\
  \EMap{f}{
    \ERet{a}
  } &=
  \ERet{f(a)}
  \\
  \EMap{f}{
    \EQuery{i}{\mathfrak{a}}
  } &=
  \EQuery{i}{
    \lambda j.\ \EMap{f}{\mathfrak{a}(j)}
  }
\end{align*}

In fact, the Escard\'o dialogues also have the structure of a monad, taking
$\ERet{-}$ as the unit:
\begin{align*}
  \EBind&: \ETree{I}{J}{A}\to \parens{A\to \ETree{I}{J}{B}} \to\ETree{I}{J}{B}
  \\
  \EBind{
    \ERet{a}
  }{f} &= f(a)
  \\
  \EBind{
    \EQuery{i}{\mathfrak{a}}
  }{f} &=
  \EQuery{i}{
    \lambda j.\
    \EBind{\mathfrak{a}(j)}{f}
  }
\end{align*}

The functor and monad laws for $\ETree{I}{J}{-}$ hold immediately by induction
on dialogue trees; in the inductive step, however, we do need function
extensionality.

\subsubsection{Running Escard\'o dialogues}
An Escard\'o dialogue $\mathfrak{a}:\ETree{I}{J}{A}$ may be executed on a
function $\alpha : I\to J$ to return a result
$\EDecode{\mathfrak{a}}{\alpha}: A$ as follows:
\begin{align*}
  \EDecode &: \ETree{I}{J}{A} \to \parens{I \to J}\to A
  \\
  \EDecode{\ERet{a}}{\alpha} &= a
  \\
  \EDecode{\EQuery{i}{\mathfrak{a}}}{\alpha} &=
  \EDecode{\mathfrak{a}(\alpha(i))}{\alpha}
\end{align*}

The following two lemmas are due to \citet{escardo:2013}.

\begin{lemma}[Naturality of execution]\label{lem:exec-naturality}
  For any $\alpha : I \to J$, the execution map $\EDecode{-}{\alpha}$ is a
  natural transformation $\ETree{I}{J}{-}\to \ArrId{\Sets}$. In other words,
  for any $f : A\to B$, the following square commutes:
  \begin{equation}
    \DiagramSquare{
      nw = \ETree{I}{J}{A},
      ne = \ETree{I}{J}{B},
      sw = A,
      se = B,
      north = \EMap{f}{-},
      west = \EDecode{-}{\alpha},
      east = \EDecode{-}{\alpha},
      south = f,
      width = 4cm,
    }
  \end{equation}
\end{lemma}

\begin{proof}
  Immediate by induction on the dialogue tree.
\end{proof}

The following lemma is needed in Section~\ref{sec:compatibility} to establish the
compatibility of the standard model of \SystemT{} with the dialogue model defined
below in Section~\ref{sec:dialogue-model}.

\begin{lemma}[Execution of Kleisli extension]\label{lem:kleisli-extension}
  Kleisli extension commutes with execution, in the sense that the following diagram commutes for every $f: A\to \ETree{I}{J}{B}$ and $\alpha : I \to J$:
  \begin{equation}
    \begin{tikzpicture}[node distance=5cm, on grid]
      \node (TA) {$\ETree{I}{J}{A}$};
      \node (TB) [right = of TA] {$\ETree{I}{J}{B}$};
      \node (A) [below = 2cm of TA] {$A$};
      \node (B) [right = of A, below = 2cm of TB] {$B$};
      \node (TB') [below = 1cm of B, between = A and B] {$\ETree{I}{J}{B}$};
      \path[->] (TA) edge node [above] {$\EBind{-}{f}$} (TB);
      \path[->] (TA) edge node [left] {$\EDecode{-}{\alpha}$} (A);
      \path[->] (TB) edge node [right] {$\EDecode{-}{\alpha}$} (B);
      \path[->] (A) edge [bend right = 10] node [sloped,below] {$f$} (TB');
      \path[->] (TB') edge [bend right = 10] node [sloped, below] {$\EDecode{-}{\alpha}$} (B);
    \end{tikzpicture}
  \end{equation}
\end{lemma}

\begin{proof}
  Immediate by induction on the dialogue tree.
\end{proof}

\subsubsection{The dialogue model of \texorpdfstring{\SystemT{}}{System T}}\label{sec:dialogue-model}

We give an effectful \emph{call-by-name} interpretation of \SystemT{}, in which
each type is interpreted as an \emph{algebra} for the dialogue monad
$\EBaire{-}$, as suggested by one of the anonymous referees. The role of
algebras here is the usual for effectful call-by-name semantics: a type is a
set equipped with the capability to make queries to an oracle.

\begin{definition}
  An algebra for the dialogue monad is a set $A$ together with a function $\Alg{A} :
  \EBaire{A}\to {A}$ satisfying the following conditions:\footnote{In fact,
  only the first condition is needed for our main result; in spite of this, it
  is more natural to speak of algebras for a monad than to speak of algebras
  for its underlying ``pointed endofunctor''. Our intention is to identify the
  most intuitive proof, not the one that makes the minimal assumptions.}
  \begin{enumerate}
    \item For all $a : A$, we have $\Alg{A}(\ERet{a}) = a$.
    \item For all $m : \EBaire{\EBaire{A}}$, we have $i(\EMap{\Alg{A}}{m}) = \Alg{A}(\EJoin{m})$ where $\EJoin{m} = \EBind{m}{\lambda x.x}$.
  \end{enumerate}
\end{definition}

\begin{remark}
  The purpose of the algebra laws is to ensure that the action is compatible with
  the monadic structure of Escard\'o trees. For instance, the first law expresses
  that a ``constant interaction'' with the oracle that always returns $a$ is
  taken to $a$; the second law ensures that interpreting a nested interaction is
  the same as interpreting each level of the interaction successively.
\end{remark}

There is a forgetful functor $\Forget$ from algebras to sets taking an algebra
$\underline{A} = (A,\Alg{A})$ to its carrier set $\Forget{\underline{A}} = A$.

\begin{construction}[Free algebras]
  The forgetful functor has a left adjoint
  $\Free$ that takes a set $A$ to the \emph{free} algebra $\Free{A}$ whose
  carrier is the set of Escard\'o trees on $A$.
  \begin{gather*}
    \Forget{\Free{A}} = \EBaire{A}
    \qquad
    \begin{aligned}
      \Alg{\Free{A}} &: \EBaire{\Forget{\Free{A}}}\to \Forget{\Free{A}}\\
      \Alg{\Free{A}}(m) &= \EBind{m}{\lambda x.x}
    \end{aligned}
  \end{gather*}
\end{construction}

\begin{construction}[Function algebras]\label{con:fun-alg}
  Given a set $A$ and an algebra $\underline{B}$, we may define the
  \emph{function algebra} $A\to\underline{B}$ whose carrier set is the function
  set $A\to\Forget{\underline{B}}$; the algebra structure is obtained by applying
  $\Alg{B}$ pointwise.
  \begin{gather*}
    \Forget{A\to\underline{B}} = A \to \Forget{\underline{B}}
    \qquad
    \begin{aligned}
      \Alg{A\to\underline{B}} &: \EBaire{\Forget{A\to\underline{B}}}\to \Forget{A\to \underline{B}}\\
      \Alg{A\to\underline{B}}(m) &= \lambda x. \Alg{\underline{B}}(\EMap{(\lambda f.f(x))}{m})
    \end{aligned}
  \end{gather*}
\end{construction}

\begin{construction}[Dependent product algebras]\label{con:dprod-alg}
  Likewise, given a set $I$ and a family of algebras
  $\underline{B}_{({i\in I})}$, we may form the dependent product algebra
  $\prod_{i\in I}{\underline{B}_i}$ in the same way.
  \begin{gather*}
    \Forget \big({\textstyle\prod_{i\in I}{\underline{B}_i}}\big)
    = {\textstyle\prod_{i\in I}\Forget{\underline{B}_i}}
    \qquad
    \begin{aligned}
      \Alg{\prod_{i\in I}{\underline{B}_i}} &:
      \EBaire{
        \Forget \big({\textstyle\prod_{i\in I}{\underline{B}_i}}\big)
      }
      \to
      \Forget \big({\textstyle\prod_{i\in I}{\underline{B}_i}}\big)
      \\
      \Alg{\prod_{i\in I}\underline{B}_i}(m) &=
      \lambda i.
      \Alg{\underline{B}_i}(\EMap{(\lambda f.f(i))}{m})
    \end{aligned}
  \end{gather*}
\end{construction}

\begin{remark}
  In both Constructions \ref{con:fun-alg} \& \ref{con:dprod-alg}, the verification of the
  algebra laws requires function extensionality --- which we already needed to
  prove the monad laws for the dialogue monad.
\end{remark}

The above is enough to define the interpretation of \SystemT{} contexts and types $\Gamma,\sigma$
into algebras $\DiaDomG{\Gamma},\DiaDomV{\sigma}$ as follows:\footnote{It is worth noting that our
interpretation of the natural numbers as the free algebra on the set $\Nat$ is somewhat bizarre from the standpoint of call-by-name semantics,
where one expects the successor nodes to encapsulate an effect. However, as
pointed out by one of the anonymous referees, one seems to disrupt the main
result by employing the call-by-name natural numbers.}
\begin{align*}
  \DiaDomG{\Gamma} &= \textstyle\prod_{x\in\Dom{\Gamma}}\DiaDomV{\Gamma(x)}\\
  \DiaDomV{\TNat} &= \Free{\Nat}\\
  \DiaDomV{\TArr{\sigma}{\tau}} &= \Forget\DiaDomV{\sigma}\to \DiaDomV{\tau}
\end{align*}

We may now define the interpretation of terms into the dialogue model:
\begin{align*}
  \DiaInterp{\Gamma}{M}{\sigma} &:
  \Forget\DiaDomG{\Gamma}\to\Forget\DiaDomV{\sigma}
  \\
  \DiaInterp{\Gamma,x:\sigma,\Delta}{x}{\sigma}\mathfrak{g} &= \mathfrak{g}(x)
  \\
  \DiaInterp{\Gamma}{\TLam{x}{M}}{\TArr{\sigma}{\tau}}\mathfrak{g} &=
  \lambda \mathfrak{a}. \DiaInterp{\Gamma,x:\sigma}{M}{\tau}\parens*{\mathfrak{g},x\mapsto \mathfrak{a}}
  \\
  \DiaInterp{\Gamma}{\TAp{\sigma}{M}{N}}{\tau}\mathfrak{g} &=
  \parens*{\DiaInterp{\Gamma}{M}{\TArr{\sigma}{\tau}}\mathfrak{g}}\parens*{
    \DiaInterp{\Gamma}{N}{\sigma}\mathfrak{g}
  }
  \\
  \DiaInterp{\Gamma}{\TZe}{\TNat}\mathfrak{g} &=
  \ERet{0}
  \\
  \DiaInterp{\Gamma}{\TSu{M}}{\TNat}\mathfrak{g} &=
  \EMap{(1+-)}{\DiaInterp{\Gamma}{M}{\TNat}\mathfrak{g}}
  \\
  \DiaInterp{\Gamma}{\TRec{N}{M_z}{x}{h}{M_s}}{\sigma}\mathfrak{g} &=
  \Alg{\DiaDomV{\sigma}}\parens{
    \EBind{\DiaInterp{\Gamma}{N}{\TNat}}{R}
  }
  \\
  \multispan2{where $
    R(n) = {
      \Rec_{\Forget\DiaDomV{\sigma}}
      \begin{pmatrix*}[l]
        n,\\
        \ERet{\DiaInterp{\Gamma}{M_z}{\sigma}\mathfrak{g}},\\
        \lambda (\mathfrak{a},\mathfrak{b}).\
        \ERet{
          \DiaInterp{\Gamma,x:\TNat,h:\sigma}{M_s}{\sigma}\parens{
            \mathfrak{g},x\mapsto \ERet{a},h\mapsto \Alg{\DiaDomV{\sigma}}{\mathfrak{b}}
          }
        }
      \end{pmatrix*}
    }
  $}
\end{align*}

\paragraph{Comparison with Escard\'o's interpretation}

Our dialogue interpretation of \SystemT{} differs from that of Escard\'o, which
fails to treat the recursor compositionally; Escard\'o interprets every type $\sigma$ as
a set $\DiaDomV{\sigma}$, and then to interpret the recursor he must define a \emph{lifting} of
the Kleisli extension to higher type by recursion on type structure:
\begin{align*}
  \langle\EBindOp\rangle_\sigma &: \EBaire{X} \to (X\to\DiaDomV{\sigma}) \to \DiaDomV{\sigma}\\
  m \mathbin{\langle\EBindOp\rangle_\TNat} f &= \EBind{m}{f}\\
  m \mathbin{\langle\EBindOp\rangle_{\TArr{\sigma}{\tau}}} f &=
  \lambda s : \DiaDomV{\sigma}.
  m \mathbin{\langle\EBindOp\rangle_\tau} \lambda x : X. f(x, s)
\end{align*}

Then recursion at $\sigma$ is interpreted by taking the lifted Kleisli
extension $\langle\EBindOp\rangle_\sigma$ of the recursor $R(-)$ defined above.
While this style of interpretation works fine, it is not compositional: the
meaning of an expression should be determined in directly from the meaning of
its subexpressions, but here we resorted to assigning a different meaning for
the recursor at each type. \citet{levy:2006} already suggests to resolve this
problem by means of the algebra translation in Section 2.1 of his monograph on
call-by-push-value, which leads to the more elegant account of recursion
presented in this paper.

\subsubsection{The generic point}\label{sec:generic-point}

In addition to the standard natural numbers (which are interpreted as
``constant'' trees, i.e.\ trees whose leaves are all labeled by the same
numeral), the dialogue model of \SystemT{} contains \emph{non-standard} natural
numbers that do not lie in the image of the interpretation function).

As a consequence of having non-standard numbers, the dialogue model also has
non-standard sequences; the most important of these sequences is called the
\emph{generic point} $\Generic:\Forget\DiaDomV{\TArr{\TNat}{\TNat}}$, a special
non-standard sequence that can be used to probe a higher-order function for
its intensional structure, providing the main ingredient for the tabulation of
higher-order functions as trees.
\begin{align*}
  \Generic &:
    \Forget\DiaDomV{\TArr{\TNat}{\TNat}}
  \\
  \Generic &=
  \lambda \mathfrak{n}.\
  \EBind{\mathfrak{n}}{
    \lambda x.\
    \EQuery{x}{\ERet{-}}
  }
\end{align*}

Intuitively, by applying the dialogue interpretation of a functional
$\IsTm{\cdot}{F}{\TArr*{\TArr{\TNat}{\TNat}}{\TNat}}$ to this generic point, we get a
dialogue tree
$\DiaInterp{\cdot}{F}{\TArr*{\TArr{\TNat}{\TNat}}{\TNat}}\parens*{\Generic}:\EBaire{\Nat}$
that is precisely the trace of $F$'s calls to its argument. Then, assuming that $F$ witnesses $\TCovers{\vec{u}}{Q}$, we can
compute the derivation of $\IndCovers{\vec{u}}{Q}$ by induction on this trace.

\begin{lemma}\label{lem:generic-diagram}
  The generic point commutes with dialogue execution in the sense that the
  following diagram commutes for all $\alpha:\Nat\to\Nat$:
  \begin{equation}
    \DiagramSquare{
      sw = \Nat,
      se = \Nat,
      nw = \Forget\DiaDomV{\TNat},
      ne = \Forget\DiaDomV{\TNat},
      north = \Generic,
      south = \alpha,
      west = \EDecode{-}{\alpha},
      east = \EDecode{-}{\alpha},
      width = 3cm,
      height = 1.5cm,
    }
  \end{equation}
\end{lemma}
\begin{proof}
  Fixing $\mathfrak{n}:\Forget\DiaDomV{\TNat}$, one obtains $\alpha\parens*{\EDecode{\mathfrak{n}}{\alpha}}
  = \EDecode{\Generic(\mathfrak{n})}{\alpha}$ immediately by induction on the dialogue
  tree $\mathfrak{n}$.
\end{proof}

\subsection{Compatibility of standard and interactive semantics}\label{sec:compatibility}

In order to use the generic point to tabulate a higher-order function
$\IsTm{\cdot}{F}{\TArr*{\TArr{\TNat}{\TNat}}{\TNat}}$ into a dialogue tree, we
will need to establish that the interactive semantics $\Angles{\ldots}$ and the
standard semantics $\Brackets{\ldots}$ are compatible at all types, in a sense
that we must define. In essence, we will be using the fact that the dialogues
in the image of the interpretation of \SystemT{} encode pure functions, even
though they coexist with ``impure'' functions.

For intuition, it is worth considering what it would mean for the two
interpretations to be compatible at the base type $\TNat$. In this case, we would want to ensure for
any term $\IsTm{\cdot}{M}{\TNat}$, that $\StdInterp{\cdot}{M}{\TNat}$ is the
same as $\EDecode{\DiaInterp{\cdot}{M}{\TNat}}{\alpha}$ for all sequences
$\alpha:\Nat\to\Nat$:
\begin{equation}\label{diagram:base-compatibility}
  \begin{tikzpicture}[node distance=5cm, on grid]
    \node (Dia) {$\Forget\DiaDomV{\TNat}$};
    \node (Std) [right = of Dia] {$\StdDomV{\TNat}$};
    \node (Syn) [above = 3cm of Dia, between = Dia and Std] {$\braces{M\mid \IsTm{\cdot}{M}{\TNat}}$};
    \path[->] (Syn) edge [above] node [sloped] {$\DiaInterp{\cdot}{-}{\TNat}$} (Dia);
    \path[->] (Syn) edge [above] node [sloped] {$\StdInterp{\cdot}{-}{\TNat}$} (Std);
    \path[->] (Dia) edge node [below] {$\EDecode{-}{\alpha}$} (Std);
  \end{tikzpicture}
\end{equation}

In order to prove the above, however, we cannot simply proceed by induction on
the syntax of \SystemT{} --- the induction hypothesis will not be strong enough
when we pass underneath a binder; the compatibility condition we have stated in in
fact a closure property of the entire language that must be established at all
higher types simultaneously, a situation that calls for \emph{logical
relations} \citep{tait:1967}.

Tait's method of logical relations is a tool to endow a property of terms at
base type with a hereditary action on all higher types. In our case, we will
extend the property of ``having compatible standard and dialogue
interpretations'' from $\TNat$ to all types $\sigma$, and then show that the
resulting structure is closed under all the constructors of \SystemT{}.

\subsubsection{A logical relation between the two models}

We begin by defining the main constituents of a logical relation between
$\StdDomV{-}$ and $\Forget\DiaDomV{-}$.
For each point $\alpha : \Nat\to\Nat$, we will define an interpretation of each
type $\sigma$ as a relation between $\StdDomV{\sigma}$ and $\Forget\DiaDomV{\sigma}$
that expresses the \emph{hereditary} compatibility of the standard
interpretation with the execution of the dialogue interpretation at $\alpha$.

.iven $s : \StdDomV{\sigma}$ and $\frk{s} : \Forget\DiaDomV{\sigma}$, we will
define the following predicate by recursion on the type $\sigma$:
\begin{mathpar}
  \boxed{
    \LR{\alpha}{\sigma}{s}{\frk{s}}
  }
  \\
  \inferrule{
    n = \EDecode{\frk{n}}{\alpha}
  }{
    \LR{\alpha}{\TNat}{n}{\frk{n}}
  }
  \and
  \inferrule{
    \forall s:\StdDomV{\sigma}, \frk{s}:\Forget\DiaDomV{\sigma}.\
    \LR{\alpha}{\sigma}{s}{\frk{s}} \implies
    \LR{\alpha}{\tau}{f(s)}{\frk{f}(\frk{s})}
  }{
    \LR{\alpha}{\TArr{\sigma}{\tau}}{f}{\frk{f}}
  }
\end{mathpar}

Likewise, we may define an analogous relation on contexts and environments as
follows, given $g:\StdDomG{\Gamma}$ and $\frk{g}:\Forget\DiaDomG{\Gamma}$:
\begin{mathpar}
  \inferrule{
    \forall x\in\Dom{\Gamma}.\
    \LR{\alpha}{\Gamma(x)}{g(x)}{\frk{g}(x)}
  }{
    \LR{\alpha}{\Gamma}{g}{\frk{g}}
  }
\end{mathpar}

\begin{theorem}[Fundamental theorem]\label{thm:fundamental-theorem}
  Fix a point $\alpha:\Nat\to\Nat$ and
  let $\IsTm{\Gamma}{M}{\sigma}$ be a term of $\SystemT{}$; then for all
  $g:\StdDomG{\Gamma}$ and $\frk{g}:\Forget\DiaDomG{\Gamma}$ such that
  $\LR{\alpha}{\Gamma}{g}{\frk{g}}$, we have
  $\LR{\alpha}{\sigma}{\StdInterp{\Gamma}{M}{\sigma}g}{\DiaInterp{\Gamma}{M}{\sigma}\frk{g}}$.
\end{theorem}

\begin{proof}
  By induction on $\IsTm{\Gamma}{M}{\sigma}$, using the monad axioms, the algebra axioms, and
  Lemmas~\ref{lem:exec-naturality}~and~\ref{lem:kleisli-extension}.
\end{proof}

\subsection{Tabulating realizable bars into Escard\'o dialogues}

Fixing a monotone subset $Q\subseteq\List{\Nat}$, our goal has been to show
that $\IndCovers{\vec{u}}{Q}$ follows from $\TCovers{\vec{u}}{Q}$, recalling
the definition of the latter:
\begin{align*}
  \parens*{\TCovers{\vec{u}}{Q}} &\triangleq
  \exists (\IsTm{\cdot}{f}{\TArr*{\TArr{\TNat}{\TNat}}{\TNat}}).\
  \forall \alpha:\Nat\to\Nat.\
  \Prepend{\vec{u}}{\Take{\TMetaAp{f}{\alpha}}{\alpha}}
  \in Q
\end{align*}

Using the results of the previous section, we are now equipped to prove an
intermediate theorem that tabulates any \SystemT{}-realizable bar into an
\emph{Escard\'o} dialogue, in the sense that $\TCovers{\vec{u}}{Q}$ implies
$\EDiaCovers{\vec{u}}{Q}$, defined below:
\begin{align*}
  \parens*{\EDiaCovers{\vec{u}}{Q}} &\triangleq
  \exists \mathfrak{n} : \EBaire{\Nat}.\
  \forall \alpha:\Nat\to\Nat.\,
  \Prepend{\vec{u}}{\Take{\EDecode{\mathfrak{n}}{\alpha}}{\alpha}}\in Q
\end{align*}

It is here that we make use of the generic point
$\Generic:\Forget\Angles{\TArr{\TNat}{\TNat}}$ in a critical way.

\begin{theorem}[Tabulation of realizable bars]\label{thm:etree-from-bar}
  If $\TCovers{\vec{u}}{Q}$ then $\EDiaCovers{\vec{u}}{Q}$.
\end{theorem}
\begin{proof}
  We abbreviate $\sigma \triangleq \TArr*{\TArr{\TNat}{\TNat}}{\TNat}$. Suppose
  $\TCovers{\vec{u}}{Q}$, i.e.\ we have some term $\IsTm{\cdot}{F}{\sigma}$
  such that $\Prepend{\vec{u}}{\Take{\TMetaAp{F}{\alpha}}{\alpha}}\in Q$ for all $\alpha:\Nat\to\Nat$.
  Interpreting $F$ into the dialogue model, we have $\mathfrak{F}\triangleq
  \Angles{F}:\Forget\Angles{\sigma}$; instantiating with the
  generic point, we therefore obtain a single Escard\'o dialogue
  $\mathfrak{F}(\Generic):\Forget\DiaDomV{\TNat}$.

  Using $\mathfrak{F}(\Generic)$ as our witness to $\EDiaCovers{\vec{u}}{Q}$,
  we fix $\alpha:\Nat\to\Nat$ to verify that $Q$ contains the list
  $\Prepend{\vec{u}}{\Take{\EDecode{\mathfrak{F}(\Generic)}{\alpha}}{\alpha}}$.
  Recalling our assumption, it would suffice to verify that
  $\TMetaAp{F}{\alpha} = \EDecode{\mathfrak{F}(\Generic)}{\alpha}$.

  By the fundamental theorem
  (Theorem~\ref{thm:fundamental-theorem}), we already have
  $\LR{\alpha}{\sigma}{\Brackets{F}}{\frk{F}}$.
  Unfolding the logical
  relation explicitly
  this amounts to the fact that $\Brackets{F}(\beta) =
  \EDecode{\mathfrak{F}(\mathfrak{b})}{\alpha}$ for all
  $\LR{\alpha}{\TArr{\TNat}{\TNat}}{\beta}{\frk{b}}$;
  we note that $\Brackets{F}(\beta) =
  \TMetaAp{F}{\beta}$. Choosing $\beta = \alpha$ and $\mathfrak{b} = \Generic$,
  it therefore remains only to show that
  $\LR{\alpha}{\TArr{\TNat}{\TNat}}{\alpha}{\Generic}$. We fix $\LR{\alpha}{\TNat}{n}{\frk{n}}$
  to show that
  $\LR{\alpha}{\TNat}{\alpha\parens{n}}{\Generic(\frk{n})}$,
  which is
  the same as to say $\alpha(n)=\EDecode{\Generic(\mathfrak{n})}{\alpha}$. But
  this is immediate by Lemma~\ref{lem:generic-diagram}.
\end{proof}

\section{Brouwerian dialogues: a normal form for Escard\'o Dialogues}

In the previous section, we showed that a \SystemT{}-realizable bar can be
tabulated into an Escard\'o dialogue: in other words,
$\EDiaCovers{\vec{u}}{Q}$ follows from $\TCovers{\vec{u}}{Q}$. In this section,
we bridge the gap between the Escard\'o dialogues and Brouwer's
$(\boldeta,\bolddigamma)$-trees, by proving that the latter express a \emph{normal
form} for the former with respect to permutation, omission, and repetition of
queries to the choice sequence.

\subsection{Ephemerality, the essence of Brouwerian dialogues}

The content of Brouwer's purported (but failed) proof of his bar thesis was to
assert that one can tabulate the evidence for barhood into a well-founded mental
construction~\citep{van-atten:2004, dummett:elements}; Escard\'o's translation
of \SystemT{} terms into dialogue trees is essentially a formalization of
Brouwer's insight.

However, Escard\'o's dialogues differ from Brouwer's mental constructions of
barhood (which are captured precisely by the judgment $\IndCovers{\vec{u}}{Q}$) in
one crucial respect: whereas Escard\'o's trees branch on an arbitrary query to the
ambient choice sequence, queries in Brouwer's mental constructions must be made
in order, i.e.\ with respect to the current moment in ideal time; as such, the
Brouwerian dialogues are \emph{ephemeral}---with each query, the head of the
ambient choice sequence is consumed and the remainder of the dialogue is
interpreted with respect to the tail of the choice sequence.

Our task, then, will be to normalize Escard\'o's dialogues into Brouwer's
ephemeral mental constructions, and then show how to massage these into a
derivation of $\IndCovers{\vec{u}}{Q}$.
Below we define the set of Brouwerian dialogues $\BTree{J}{A}$ (coding
functionals $\parens{\Nat\to J}\to A$) as the least set closed under the rules below:
\begin{mathpar}
  \boxed{
    \inferrule{
      J,A:\Sets
    }{
      \BTree{J}{A}:\Sets
    }
  }
  \and
  \inferrule[\RuleBRet]{
    a:A
  }{
    \BRet{a}:\BTree{J}{A}
  }
  \and
  \inferrule[\RuleBQuery]{
    \mathbf{a}:\parens{J\to\BTree{J}{A}}
  }{
    \BQuery{\mathbf{a}}: \BTree{J}{A}
  }
\end{mathpar}

Just as we showed how to execute an Escard\'o dialogue against a choice
sequence, we can do the same for the
Brouwerian (ephemeral) version. For $\mathbf{a}: \BTree{J}{A}$ and
$\alpha : \Nat\to  J$, we define $\BDecode{\mathbf{a}}{\alpha}:A$ by
recursion on $\mathbf{a}$ as follows:
\begin{align*}
  \BDecode &: \BTree{J}{A}\to \parens{\Nat\to J}\to A \\
  \BDecode{\BRet{a}}{\alpha} &= a\\
  \BDecode{\BQuery{\mathbf{a}}}{\alpha} &=
  \BDecode{\mathbf{a}(\Head{\alpha})}{\Tail{\alpha}}
\end{align*}

\subsubsection{Barhood relative to Brouwerian dialogues}

Just as we have defined $\EDiaCovers{\vec{u}}{Q}$ as a notion of barhood
witnessed by an Escard\'o dialogue, we now do the same for Brouwerian
dialogues, recalling the definitions of barhood we have given so far:
\begin{align*}
  \parens*{\TCovers{\vec{u}}{Q}} &\triangleq
  \exists (\IsTm{\cdot}{f}{\TArr*{\TArr{\TNat}{\TNat}}{\TNat}}).\
  \forall \alpha:\Nat\to\Nat.\
  \Prepend{\vec{u}}{\Take{\TMetaAp{f}{\alpha}}{\alpha}}
  \in Q
  \\
  \parens*{\EDiaCovers{\vec{u}}{Q}} &\triangleq
  \exists \mathfrak{n}:\EBaire{\Nat}.\
  \forall \alpha:\Nat\to\Nat.\,
  \Prepend{\vec{u}}{\Take{\EDecode{\mathfrak{n}}{\alpha}}{\alpha}}\in Q
  \\
  {\color{Red}\parens*{\BDiaCovers{\vec{u}}{Q}}} &\triangleq
  {\color{Red}
  \exists \mathbf{n}:\BTree{\Nat}{\Nat}.\
  \forall \alpha:\Nat\to\Nat.\,
  \Prepend{\vec{u}}{\Take{\BDecode{\mathbf{n}}{\alpha}}{\alpha}}\in Q
  }
\end{align*}

The Brouwerian dialogues that we have defined are nothing but representations
of proofs of $\IndCovers{\vec{u}}{Q}$, as the following lemma shows.

\begin{lemma}\label{lem:diag-cover-to-ind-cover}
  For any monotone subset $\IsSubsetEq{Q}{\List{\Nat}}$ and basic observation
  $\vec{u}:\List{\Nat}$, if $\BDiaCovers{\vec{u}}{Q}$ then $\IndCovers{\vec{u}}{Q}$.
\end{lemma}
\begin{proof}
  Fixing $Q$ and $\vec{u}$ such that $\BDiaCovers{\vec{u}}{Q}$, we therefore have some
  Brouwerian dialogue $\mathbf{n}:\BTree{\Nat}{\Nat}$ such that for all sequences
  $\alpha$, we have $\Prepend{\vec{u}}{\Take{\BDecode{\mathbf{n}}{\alpha}}{\alpha}}\in Q$. We
  proceed by induction on $\mathbf{n}$ to prove $\IndCovers{\vec{u}}{Q}$, using the fact
  that $Q$ is monotone.
\end{proof}

\subsection{Normalizing Escard\'o dialogues into Brouwerian dialogues}

We wish to normalize elements of $\ETree{\Nat}{J}{A}$ (Escard\'o dialogues)
into elements of $\BTree{J}{A}$ (Brouwerian dialogues) in a way that commutes
with execution. To do so, we design an inductive/proof-theoretic
characterization of the normalizable Escard\'o dialogues, and then show that
all Escard\'o dialogues can be coded as such.  This yields a constructive and
structurally recursive normalization algorithm.  To this end, we define below
two mutually inductive forms of judgment, whose rules are given in
Fig.~\ref{fig:normalization}:

\begin{enumerate}
  \item $\CanNorm{\vec{u}}{\mathfrak{a}}{\mathbf{a}}$, presupposing
    $\vec{u}:\List{J}$ and $\mathfrak{a}:\ETree{\Nat}{J}{A}$, and
    guaranteeing $\mathbf{a}:\BTree{J}{A}$, means that the Escard\'o dialogue
    $\mathfrak{a}$ normalizes to the Brouwerian dialogue $\mathbf{a}$.

  \item $\CanNormQ{\vec{u}}{i}{\mathfrak{a}}{\vec{v}}{\mathbf{a}}$ presupposes
    $\vec{u},\vec{v}:\List{J}$, $i: \Nat$ and
    $\mathfrak{a}:J\to\ETree{\Nat}{J}{A}$, and guarantees
    $\mathbf{a}:\BTree{J}{A}$.
\end{enumerate}

For the sake of intuition, we provide an implementation of the algorithm in
Standard~ML in Fig.~\ref{fig:sml-normalization}; the mathematical version
given in this section can be seen as a specification and termination proof for
the Standard~ML program.

\subsubsection{Explanation of algorithm}

Normalization occurs with respect to a node $\vec{u}$ that is extended every
time a $\bolddigamma$ node is inserted; $\vec{u}$ should be thought of as the
``current known prefix of the choice sequence''. To normalize an Escard\'o
dialogue $\mathfrak{a}$ over a node $\vec{u}$, one proceeds by case: if
$\mathfrak{a} = \ERet{a}$ then it is already normal; on the other hand, if
$\mathfrak{a}=\EQuery{i}{\mathfrak{b}}$, then it is necessary to find out what
the $i$th element of the choice sequence is.

If $i<\Len{\vec{u}}$, then we have already ``cached'' this query in
$\vec{u}_i$, and we therefore proceed by normalizing $\mathfrak{b}(\vec{u}_i)$;
on the other hand, if $i\geq\Len{\vec{u}}$, we must continue to insert
$\bolddigamma$-queries until we have information about the $i$th element of the
choice sequence. Locating ``cached'' information and inserting new queries is
the role of the $\CanNormQ{\vec{u}}{i}{\mathfrak{a}}{\vec{v}}{\mathbf{a}}$
judgment.

\begin{figure*}
  \begin{mathpar}
    \DeclJdg{\CanNorm{\vec{u}}{\mathfrak{a}}{\mathbf{a}}}{
      \vec{u}:\List{J},
      \mathfrak{a}:\ETree{\Nat}{J}{A},
      \mathbf{a}:\BTree{J}{A}
    }
    \and
    \DeclJdg{\CanNormQ{\vec{u}}{i}{\mathfrak{a}}{\vec{v}}{\mathbf{a}}}{
      \vec{u},\vec{v}:\List{J},
      i:\Nat,
      \mathfrak{a}:J\to\ETree{\Nat}{J}{A},
      \mathbf{a}:\BTree{J}{A}
    }
    \\
    \inferrule{
    }{
      \CanNorm{\vec{u}}{\ERet{a}}{\BRet{a}}
    }
    \and
    \inferrule{
      \CanNormQ{\vec{u}}{i}{\mathfrak{a}}{\vec{u}}{\mathbf{a}}
    }{
      \CanNorm{\vec{u}}{\EQuery{i}{\mathfrak{a}}}{\mathbf{a}}
    }
    \and
    \inferrule{
      \CanNorm{\vec{u}}{\mathfrak{a}(j)}{\mathbf{a}}
    }{
      \CanNormQ{\vec{u}}{0}{\mathfrak{a}}{\Cons{j}{\vec{v}}}{\mathbf{a}}
    }
    \and
    \inferrule{
      \CanNormQ{\vec{u}}{i}{\mathfrak{a}}{\vec{v}}{\mathbf{a}}
    }{
      \CanNormQ{\vec{u}}{i + 1}{\mathfrak{a}}{\Cons{j}{\vec{v}}}{\mathbf{a}}
    }
    \and
    \inferrule{
      \forall j: J.\
      \CanNorm{\Snoc{\vec{u}}{j}}{\mathfrak{a}(j)}{\mathbf{a}(j)}
    }{
      \CanNormQ{\vec{u}}{0}{\mathfrak{a}}{\Nil}{\BQuery{\mathbf{a}}}
    }
    \and
    \inferrule{
      \forall j: J.\
      \CanNormQ{\Snoc{\vec{u}}{j}}{i}{\mathfrak{a}}{\Nil}{\mathbf{a}(j)}
    }{
      \CanNormQ{\vec{u}}{i + 1}{\mathfrak{a}}{\Nil}{\BQuery{\mathbf{a}}}
    }
  \end{mathpar}
  \caption{Rules for dialogue normalization.}\label{fig:normalization}
\end{figure*}

\begin{figure*}
  \begin{code}
    datatype ('i, 'j, 'a) etree =
       RET of 'a
     | QUERY of 'i * ('j $\to$ ('i, 'j, 'a) etree)

    \-
    datatype ('j, 'a) btree =
       SPIT of 'a
     | BITE of ('j $\to$ ('j, 'a) btree)

    \-
    fun norm us (RET a) = SPIT a
      | norm us (QUERY (i, k)) = normQ us i k us

    \-
    and normQ us Z k [] = BITE (fn u $\Rightarrow$ norm (us @ [u]) (k u))
      | normQ us Z k (v $::$ \_) = norm us k v
      | normQ us (S i) (\_ $::$ vs) = normQ us i k vs
      | normQ \_ \_ \_ \_ = raise Subscript
  \end{code}
  \caption{A Standard ML implementation of the dialogue normalization algorithm whose graph is specified in Fig.~\ref{fig:normalization}.}\label{fig:sml-normalization}
\end{figure*}

\begin{lemma}\label{lem:norm-complete}
  The inductive characterization of normalization
  $\CanNorm{\vec{u}}{\mathfrak{a}}{\mathbf{a}}$ in Fig.~\ref{fig:normalization} is functional, i.e.\ for any
  $\vec{u}:\List{J}$ and $\mathfrak{a}:\ETree{\Nat}{J}{A}$ there is a
  unique $\mathbf{a}:\BTree{J}{A}$ such that $\CanNorm{\vec{u}}{\mathfrak{a}}{\mathbf{a}}$.
\end{lemma}

\begin{proof}
  In fact, we must simultaneously establish the functionality of
  $\CanNorm{\vec{u}}{\mathfrak{a}}{\mathbf{a}}$ and
  $\CanNormQ{\vec{u}}{i}{\mathfrak{b}}{\vec{v}}{\mathbf{b}}$. The existence
  part of functionality follows by simultaneous induction on $\mathfrak{a}$,
  $\vec{v}$, and $i$; the uniqueness part follows immediately by inspection of
  the generating rules.
\end{proof}

\begin{corollary}[Normalization function]\label{cor:norm}
  We have a structurally recursive function $\Norm{\vec{u}}{\mathfrak{a}}$ such
  that for all $\vec{u}:\List{J}$ and $\mathfrak{a}:\ETree{\Nat}{Y}{Z}$,
  $\CanNorm{\vec{u}}{\mathfrak{a}}{\Norm{\vec{u}}{\mathfrak{a}}}$.
\end{corollary}

\begin{lemma}[Coherence of normalization and execution]\label{lem:dialogue-normalization-compatibility}
  Execution of Brouwerian dialogues is compatible with normalization, as defined in
  Corollary~\ref{cor:norm}; to be precise, the following diagram commutes for
  all $\vec{u}:\List{J}$ and $\alpha:\Nat\to J$.
  \begin{equation}\label{diagram:normalization-compatibility}
    \begin{tikzpicture}[node distance=5cm,on grid]
      \node (E) {$\ETree{\Nat}{J}{A}$};
      \node (B) [right = of E] {$\BTree{J}{A}$};
      \node (A) [below = 2cm of E, between = E and B] {$A$};
      \path[->] (E) edge node [above] {$\Norm{\vec{u}}{-}$} (B);
      \path[->] (E) edge [below] node [sloped] {$\EDecode{-}{\Prepend{\vec{u}}{\alpha}}$} (A);
      \path[->] (B) edge [below] node [sloped] {$\BDecode{-}{\alpha}$} (A);
    \end{tikzpicture}
  \end{equation}
  When $\vec{u}=\Nil$, this becomes the statement that $\EDecode{\mathfrak{a}}{\alpha}=\BDecode{\Norm{\Nil}{\mathfrak{a}}}{\alpha}$ for any dialogue $\mathfrak{a}$.
\end{lemma}
\begin{proof}
  By induction on the graph of the normalization function.
\end{proof}

\section{Validity of the bar thesis for \texorpdfstring{$\mathbf{T}$}{T}-realizable bars}\label{sec:bar-thesis-for-t}

We obtain the main theorem of this paper as a corollary to the constructions of
the past several sections.

The bar thesis for \SystemT{} states that for any monotone set
$\IsSubsetEq{Q}{\List{\Nat}}$ of nodes, the
inductive definition of barhood is complete in the sense that we can
conclude $\IndCovers{\Nil}{Q}$ from $\TCovers{\Nil}{Q}$.

\begin{theorem}[The bar thesis for realizable bars]
  Brouwer's bar thesis holds for bars that are realizable in \SystemT: for any
  monotone subset $\IsSubsetEq{Q}{\List{\Nat}}$, if $\TCovers{\Nil}{Q}$ then
  $\IndCovers{\Nil}{Q}$.
\end{theorem}

\begin{proof}
  Supposing $\TCovers{\Nil}{Q}$, we need to construct $\IndCovers{\Nil}{Q}$; by
  Lemma~\ref{lem:diag-cover-to-ind-cover}, it is enough to show that
  $\BDiaCovers{\Nil}{Q}$, i.e.\ exhibit some $\mathbf{n}:\BTree{\Nat}{\Nat}$
  such that for all $\alpha:\Nat\to\Nat$, we have
  $\Take{\BDecode{\mathbf{n}}{\alpha}}{\alpha}\in Q$.
  On the other hand, by Theorem~\ref{thm:etree-from-bar} we already have
  $\EDiaCovers{\Nil}{Q}$, from which we obtain some
  $\mathfrak{n}:\EBaire{\Nat}$ such that for all $\alpha:\Nat\to\Nat$, we
  have $\Take{\EDecode{\mathfrak{n}}{\alpha}}{\alpha}\in Q$. Therefore,
  choosing $\mathbf{n} = \Norm{\Nil}{\mathfrak{n}}$, all that remains is check
  that $\Take{\BDecode{\Norm{\Nil}{\mathfrak{n}}}{\alpha}}{\alpha} =
  \Take{\EDecode{\mathfrak{n}}{\alpha}}{\alpha}$. But
  $\BDecode{\Norm{\Nil}{\mathfrak{n}}}{\alpha} =
  \EDecode{\mathfrak{n}}{\alpha}$ by
  Lemma~\ref{lem:dialogue-normalization-compatibility}.
\end{proof}

\section{Discussion and related work}

\subsection{The Bar Induction principle}

The essence of Brouwer's bar thesis lies in the tabulation of a proof of
$\Covers{\vec{u}}{Q}$ into a proof of $\IndCovers{\vec{u}}{Q}$. In
constructive metamathematics, however, an equivalent presentation of the Bar
Thesis as an induction principle is usually considered:

\begin{proposition}[Monotone Bar Induction]\label{prop:bi-mono}

  For any subset $R\in\List{\Nat}$ and a monotone bar $\Covers{\Nil}{Q}$, we
  can conclude $\Nil\in R$ from the following conditions:
  \begin{enumerate}
    \item \emph{Base case.} $R$ contains the bar, i.e.\ $Q\subseteq R$.
    \item \emph{Inductive step.} $R$ is inductive, i.e.\ $\vec{u}\in R$ follows from $\forall x.\ \Snoc{\vec{u}}{x}\in R$.
  \end{enumerate}

\end{proposition}

To prove the equivalence of Proposition~\ref{prop:bi-mono} with our formulation
of Brouwer's bar thesis on monotone bars (Proposition~\ref{prop:completeness}),
we begin by verifying that it suffices to consider unqualified barhood of the
root node $\Nil$. In this section we write $\vec{v}\preceq\vec{u}$ to mean that
$\vec{u}=\Prepend{\vec{v}}{\vec{w}}$ for some $\vec{w}$; this notation is the
reverse of Dummett's~(\citeyear{dummett:elements}), but we find it more concrete.

\begin{lemma}[It suffices to consider the root]\label{lem:root-suffices}
  If $\Covers{\Nil}{Q}\implies\IndCovers{\Nil}{Q}$ for all monotone $Q$, then $\Covers{\vec{u}}{Q}\implies\IndCovers{\vec{u}}{Q}$ for all $\vec{u}$ and monotone $Q$.
\end{lemma}

\begin{proof}
  Suppose that $\Covers{\Nil}{Q}\implies\IndCovers{\Nil}{Q}$ for all monotone
  $Q$. Fixing monotone $Q'$ such that $\Covers{\vec{u}}{Q'}$, we need to show
  that $\IndCovers{\vec{u}}{Q'}$. Letting $Q =
  \braces{\vec{v}\mid\vec{u}\preceq\vec{v}\implies\vec{v}\in Q'}$, we observe
  that $\Covers{\Nil}{Q}$ immediately, whence $\IndCovers{\Nil}{Q}$ follows by
  assumption. To see that $\IndCovers{\vec{u}}{Q'}$ follows from $\IndCovers{\Nil}{Q}$,
  we observe that the latter can be obtained immediately as a subtree of the former.
\end{proof}

\begin{theorem}
  Monotone Bar Induction (Proposition~\ref{prop:bi-mono}) is equivalent to
  Brouwer's Thesis on monotone bars (Proposition~\ref{prop:completeness}).
\end{theorem}

\begin{proof}
  ($\Rightarrow$) Suppose that Proposition~\ref{prop:bi-mono} holds; by
  Lemma~\ref{lem:root-suffices}, it suffices to show that $\IndCovers{\Nil}{Q}$
  follows from $\Covers{\Nil}{Q}$. We choose a \emph{motive of bar induction}
  based on our goal, namely $R\triangleq
  \braces*{\vec{u}\mid\IndCovers{\vec{u}}{Q}}$.

  \begin{enumerate}
    \item \emph{Base case.} Fixing $\vec{u}\in Q$, we need to see that $\IndCovers{\vec{u}}{Q}$; this is just the \boldeta-inference.
    \item \emph{Inductive step.} Suppose that $\IndCovers{\Snoc{\vec{u}}{x}}{Q}$ for all $x$; we conclude $\IndCovers{\vec{u}}{Q}$ by the \bolddigamma-inference.
  \end{enumerate}

  ($\Leftarrow$) Suppose that Proposition~\ref{prop:completeness} holds. Fix a
  monotone bar $\Covers{\Nil}{Q}$ and a subset $R\supseteq{Q}$ such that
  $\vec{u}\in R$ follows from $\forall x.\ \Snoc{\vec{u}}{x}\in R$ for all
  $\vec{u}$, to show that $\Nil\in\vec{R}$. By
  Proposition~\ref{prop:completeness}, we have $\IndCovers{\Nil}{Q}$; we
  proceed by induction, generalizing to conclude $\vec{v}\in R$ from
  $\IndCovers{\vec{v}}{Q}$ for all $\vec{v}$:
  \begin{enumerate}
    \item \emph{\boldeta-inference.} If $\vec{v}\in Q$, then $\Nil\in R$ because $Q\subseteq R$.
    \item \emph{\bolddigamma-inference.} Suppose that $\IndCovers{\Snoc{\vec{v}}{x}}{Q}$ for all $x$; by induction, we have $\Snoc{\vec{v}}{x}\in R$ for all $x$, whence by assumption we have $\vec{v}\in R$.\qedhere
  \end{enumerate}
\end{proof}

\subsection{Monotonicity, decidability, and continuity}

We have focused exclusively on \emph{monotone} subsets $Q\subseteq\List{\Nat}$,
i.e.\ collections of basic observations that are closed under extension;
intuitively, such a bar is one that cannot be ``escaped''. However, another
variant of the bar thesis that concerns \emph{decidable} subsets has received
more attention in the literature.
A subset $Q\subseteq\List{\Nat}$ is decidable when $\forall\vec{u}.\ \vec{u}\in
Q \lor \vec{u}\not\in Q$ holds; of course, this is classically true for every
subset $Q$, but in constructive logic, it doesn't necessarily hold for all $Q$.

The bar induction principle for decidable bars is weaker than monotone bar
induction (Proposition~\ref{prop:bi-mono}); in essence, this is because any
decidable bar can be completed into a monotone bar~\citep{dummett:elements}.

\begin{lemma}\label{lem:bar-free-mono}
  Any decidable bar $\Covers{\Nil}{Q}$ can be freely completed into a monotone
  decidable bar $\Covers{\Nil}{Q^\natural}$ with $Q\subseteq Q^\natural$.
\end{lemma}

\begin{proof}
  The inclusion of the \emph{discrete} poset of lists of natural
  numbers into the poset of lists with the approximation ordering
  $\Prepend{\vec{u}}{\vec{v}}\succeq\vec{u}$ induces a contravariant reindexing
  from monotone subsets into general subsets. This reindexing has both left and
  right adjoints, by Kan extension. In particular, the left adjoint takes a
  general subset $Q$ and transforms it into a monotone subset
  $Q^\natural\triangleq\braces*{\vec{u}\mid\exists\vec{v}\preceq\vec{u}.\ \vec{v}\in Q}$.

  First we observe that the transformation described above preserves the
  decidability of the bar: to decide
  $\vec{u}\in Q^\natural$, one simply decides $\vec{v}\in Q$ for as many
  prefixes $\vec{v}\preceq\vec{u}$ as necessary.
  Second, we check that $\Covers{\Nil}{Q^\natural}$ assuming
  $\Covers{\Nil}{Q}$. Fixing $\alpha:\Nat\to\Nat$, we need some $k:\Nat$ such
  that $\Take{k}{\alpha}\in Q^\natural$. But $Q$ is a bar, so we already have
  $k$ such that $\Take{k}{\alpha}\in Q$, from which we immediately obtain
  $\Take{k}{\alpha}\in Q^\natural$.
\end{proof}

\begin{theorem}
  Decidable bar induction follows from monotone bar induction.
\end{theorem}
\begin{proof}
  We fix a \emph{decidable} bar $\Covers{\Nil}{Q}$ and an inductive superset
  $R\supseteq{}Q$, to show that $\Nil\in{}R$. By Lemma~\ref{lem:bar-free-mono},
  we have a monotone decidable bar $Q^\natural$; following
  \citet{dummett:elements}, we apply monotone bar induction on
  $Q^\natural$, choosing the motive $S\triangleq R\cup Q^\natural$. We clearly have
  the base case $Q^\natural\subseteq S$; to see that $S$ is inductive, we assume
  that $\forall x.\ \Snoc{\vec{u}}{x}\in S$, to show that $\vec{u}\in S$.
  Because $Q^\natural$ is decidable, we can proceed by case:
  \begin{enumerate}

    \item If $\vec{u}\in Q^\natural$, then we immediately have $\vec{u}\in S$.

    \item If $\vec{u}\not\in Q^\natural$, then our only option is to check that
      $\vec{u}\in R$. Using the inductiveness of $R$, it suffices to verify
      that $\Snoc{\vec{u}}{x}\in R$ for all $x$. By assumption, we have
      $\Snoc{\vec{u}}{x}\in S$; we proceed by case on whether $\Snoc{\vec{u}}{x}\in Q^\natural$:
      \begin{enumerate}

        \item Suppose $\Snoc{\vec{u}}{x}\in Q^\natural$; because $\vec{u}\not\in
          Q^\natural$, we know that no prefix of $\vec{u}$ is in $Q$, so therefore
          we must have $\Snoc{\vec{u}}{x}\in Q$. Because $Q$ is a subset of
          $R$, we are done.

        \item Suppose $\Snoc{\vec{u}}{x}\not\in Q^\natural$; because we have assumed
          $\Snoc{\vec{u}}{x}\in S$, we therefore already have
          $\Snoc{\vec{u}}{x}\in R$.\qedhere

      \end{enumerate}
  \end{enumerate}
\end{proof}

Under a further continuity principle, which \citet{escardo:2013}
has shown to be validated for \SystemT*{}-definable functionals, monotone and
decidable bar induction are in fact
equivalent~\citep{dummett:elements,troelstra-vandalen:1988}.

\subsection{Schwichtenberg's closure theorem}

In 1979, Schwichtenberg proved an even stronger result than what we
have proved here, namely that for any closed \SystemT{} term that
codes a bar of type $0$ and $1$, the \emph{bar recursor} can already
be defined in \SystemT*{}~\citep{schwichtenberg:1979}. In more recent
work, Oliva and Steila give an elegant and \emph{direct} proof of
Schwichtenberg's result~\citep{oliva-steila:2016}.
It should be possible to replicate this result in our setting by using
Church encodings of dialogues, which Escard\'o has used to exhibit the modulus
of continuity of a \SystemT*-definable functional as a program in
\SystemT{}~\citep{escardo:2013}.

The results of Schwichtenberg, Oliva and Steila, Escard\'o, Xu, and ourselves
are part of long tradition of comparing formal definability with continuity and
bar induction principles, initiated by \citet{kreisel-troelstra:1970} in a
pioneering tour de force that clarified the closure conditions enjoyed by
neighborhood functions (functions that correspond to Brouwer's inductive bars).

\subsection{Bar recursion and continuous moduli}

Whereas bar induction is a principle of logic (a method to obtain proofs), bar
\emph{recursion} is an principle of mathematical construction (a method to
define functions) proposed by \citet{spector:1962}.  In results
obtained subsequent to those described in our own paper, \citet{fujiwara-kawai:2019}
show that bar induction for decidable bars is equivalent to the existence of a
bar recursor for every functional of type $(\Nat\to\Nat)\to\Nat$ that has a
\emph{continuous} modulus of continuity. This result
is closely related to that of \citet{schwichtenberg:1979}, as
well as those of Escard\'o and ourselves, considering Escard\'o's proof that
every definable such functional in \SystemT{} has a definable modulus of
continuity --- and hence a continuous modulus of
continuity~\citep{escardo:2013}.

In other recent work, also obtained subsequent to our results,\footnote{Our
proof was completed and announced in Spring of 2016.} Xu has extended
Escard\'o's methods to achieve a monadic translation of \SystemT{} into itself
parameterized in a monad; Xu's result includes an elegant
algorithm to obtain both moduli of continuity and bar recursors of \SystemT{}
programs~\citep{xu:2020:fscd,xu:2020:lmcs}.

\subsection{Martin-L\"of's analysis of \texorpdfstring{$\Pi^1_1$}{Pi1/1} sentences}

In his doctoral dissertation, Per Martin-L\"of pursued an explanation of constructive
meaning in which $\Pi^1_1$ sentences ending in a decidable predicate are
\emph{defined} to be canonically verified by $(\boldeta,\bolddigamma)$-trees as
envisioned by Brouwer~\citep{martin-lof:1970}. A $\Pi^1_1$ sentence is one
of the form $\forall \overline{\alpha_i : S_i}.\ \exists \overline{\beta_j :
T_j}.\ P\parens[\big]{\overline{\alpha_i}\ldots\overline{\beta_i}}$, in which
$S_i,T_i$ are types of order at most $1$; the statement that some decidable
subset $Q$ is a bar is an example of such a sentence, but
Martin-L\"of's insight was that one might as well consider \emph{all} such
sentences as being about bars.

In contrast to the Brouwer--Heyting--Kolmogorov (BHK) interpretation of
Intuitionistic logic, Martin-L\"of did not at the time consider every
well-formed sentence generated from
$\forall,\exists,\lor,\land,\supset,\lnot,P$ to be constructively meaningful;
rather, he argued that different kinds of sentences needed to be analyzed
individually and interpreted, in order to be endowed with constructive meaning.
This hesitant and open-ended perspective on constructive meaning is in fact
much closer to Brouwer's views than the perhaps infelicitously named BHK
interpretation --- an interpretation that, in generality, fails to verify
Brouwer's thesis on bars (which might as well be called ``Brouwer's thesis on
$\Pi^1_1$ sentences'').

Exhibiting a theory of meaning that is at once compatible with both the BHK
interpretation \emph{and} Brouwer's interpretation of $\Pi^1_1$ sentences is a
difficult matter, and leads inevitably to \emph{forcing} --- both
sheaf-theoretic and effectful. The subject of this paper is to explain one
possible BHK-compatible theory of meaning that also validates Brouwer's
thesis, but we must admit that realizability in \SystemT{} is too formalistic a
constraint to be realistic.

\subsection{Operational forcing}

A forerunner of Escard\'o's effectful forcing can be found in Coquand and
Jaber's \emph{A computational interpretation of forcing in Type
Theory}~(\citeyear{coquand-jaber:2012}); unlike effectful forcing, which is a
purely denotational construction, Coquand and Jaber give an \emph{operational}
interpretation of \SystemT{} extended by an oracle
$\mathsf{f}:\TArr{\TNat}{\TKwd{bool}}$, and use it to show that the definable
functionals on sequences of booleans are \emph{uniformly} continuous.

Using similar ideas, \citet{coquand-mannaa:2016} have constructed a countermodel to
Markov's principle for dependent type theory
from an operational perspective, adding a new constant
$\mathsf{f:\TArr{\TNat}{\TKwd{bool}}}$ to the language. In order to specify the
equational behavior of the generic point $\mathsf{f}$, the judgments of the
type theory are indexed in a \emph{forcing condition} $p$ that specifies a
finite subgraph of a sequence of booleans. Then, two conditions are imposed:
\begin{enumerate}

  \item The judgments of type theory must be \emph{monotone} with respect to
    approximations of forcing conditions.

  \item The judgments of type theory must be \emph{local} with respect to
    coverings of forcing conditions: roughly, if the judgment $\mathcal{J}$
    holds at each leaf of an $(\boldeta,\boldbeta)$ tree rooted at $p$, then
    $\mathcal{J}$ holds at $p$.

\end{enumerate}

The locality condition above has the result of rendering the semantic evidence
of the \emph{judgments} of the forcing extension into
$(\boldeta,\boldbeta)$ trees; this is to be contrasted with the denotational
style of effectful forcing, in which the terms themselves are interpreted as
$(\boldeta,\boldbeta)$ trees (potentially at higher type). This is concordant
with the usual differences between operational and denotational approaches to semantics.

\subsubsection{Forcing in Nuprl}

Nuprl is one of the oldest implementations of dependent type theory, inspired
by the ideas of Bishop, De Bruijn, Martin-L\"of, and
Scott~\citep{constable:1986}. In contrast to the algebraic tradition of type theory, in
which the sound and complete interpretation of a formal language into a
suitable class of models is emphasized, Nuprl is based on the idea of
considering a \emph{single} intended model, whose properties are meant to
approximate the limitations of human mental construction.

Over the years, the designers of Nuprl have increasingly emphasized the
validity of non-constructive but computationally justified principles in their
evolving computational model, including type-theoretic variants of Markov's
principle, bar induction, and Brouwer's continuity
theorems~\citep{constable:2014,rahli-bickford:2016,rahli-bickford-constable:2017}.
The validity of Markov's principle and bar induction in Nuprl followed roughly
from the fact that Nuprl's computational semantics had been developed in a
classical metatheory, in which these principles already hold. At a high level,
the classical metatheory of the operational model allowed termination
properties of programs to be established using classical logic.

Finding this state of affairs unsatisfactory, researchers working on Nuprl have
recently developed a new constructive semantics based on
forcing~\citep{bickford-cohen-constable-rahli:2018}, in a manner that
generalizes the forcing conditions of \citet{coquand-mannaa:2016} to include
finite data about an unbounded collection of choice sequences.

The resulting Nuprl system can be seen as a synthesis of Brouwer's purely
denotational understanding of constructive meaning with the operational
approach to constructivity pioneered by Martin-L\"of in his
influential report, \emph{Constructive Mathematics and Computer
Programming}~(\citeyear{martin-lof:1979}).

\subsubsection{Relation to sheaf models}

The operational account of forcing extensions of type theory described above is
tantalizing close to the notion of a sheaf model over a topological space (such
as the space of sequences of naturals or booleans). While it was perhaps
unclear at the time, due to the highly unfolded character of the constructions,
this connection can be rationalized in the following way: the extended
operational model of Coquand and Mannaa is, at a high level, just the
\emph{standard} operational model of dependent type theory carried out
internally to the topos of sheaves on Cantor space. A similar statement should
be true of the operational model of Nuprl's forcing
extension~\citep{bickford-cohen-constable-rahli:2018}.

More recently, \citet{sterling-harper:2018:lics} have developed an analogous
construction of a different forcing extension of type theory for \emph{guarded
recursion}, utilizing the internal language of another sheaf topos to express
the operational semantics and relational interpretation of the judgments of
type theory.

\subsection{Denotational forcing}

The main ideas underlying our proof of Brouwer's bar thesis for \SystemT*-realizable
bars were invented by Escard\'o in his paper \emph{Continuity of G\"odel's
System T definable functionals via effectful forcing}~\citep{escardo:2013}.
Escard\'o used the effectful forcing technique to prove a related property of
\SystemT*-definable functionals $F:\parens{\Nat\to\Nat}\to\Nat$,
namely continuity:
\begin{equation*}
  \forall \alpha.\
  \exists k.\
  \forall \beta.\
  \Take{k}{\alpha} = \Take{k}{\beta}
  \implies
  F(\alpha) = F(\beta)
\end{equation*}
The proof of this fact, being constructive, contains an \emph{algorithm} to
compute this modulus of continuity $k$ from each sequence $\alpha$.

Relative to Escard\'o, we have shown that the same method can be used to
validate a restricted version of Brouwer's bar thesis, contributing a novel
normalization theorem for Escard\'o's non-canonical
$(\boldeta,\boldbeta)$-trees into Brouwer's canonical $(\boldeta,\bolddigamma)$
trees. Unlike Escard\'o, who developed effectful forcing for a version of
\SystemT{} extended by an oracle, we have avoided adding the oracle to the
syntax, following a suggestion from Vincent Rahli; finally, following a
suggestion from one of the anonymous referees, we have achieved a more
compositional dialogue interpretation.

\subsubsection{Syntactic forcing translations}

Recently a number of researchers have investigated forcing extensions of the
Coq proof assistant and its calculus of inductive constructions, further
clarifying the connection between forcing, computational
effects, and parametricity~\citep{jaber-tabareau-sozeau:2012,
jaber-lewertowski-pedrot-sozeau-tabareau:2016,pedrot-tabareau:2017,pedrot-tabareau:2020}.
Coq's forcing apparatus differs from that of Nuprl in a couple important ways.

First, forcing over an arbitrary preorder is supported, whereas Nuprl has been
hard-coded for a specific forcing extension. On the other hand, Nuprl's forcing
is \emph{localized} in the sense that one may amalgamate objects defined on a
bar, in essence a sheaf condition; the local character of the Nuprl semantics
is very important, as it is the main ingredient to justifying the bar induction
principle in the forcing extension.

Second, the consistency of Coq's forcing extension is justified by purely
syntactical means through a simple translation into the base calculus. In
contrast, Nuprl's forcing extension is justified by means of a completely new
version of the partial equivalence relation semantics in which monotonicity and
local character conditions are fully unfolded \`a la
\citet{coquand-mannaa:2016} rather than treated abstractly \`a la
\citet{sterling-harper:2018:lics}, an enormous undertaking that took multiple
years.

\subsubsection{Sheaf-theoretic forcing}

An arguably more direct approach to obtain the compatiblity of type theories
with Brouwerian principles such as bar induction and continuity is to consider
the petit topos of sheaves on a suitable space (such as Baire space or Cantor
space), or the gros topos of sheaves over a suitable category of test
spaces~\citep{fourman:1984}.

In his doctoral thesis, \citet{xu:2015} considers sheaf models of dependent type
theory validating uniform continuity and the fan theorem. There
is an apparent difficulty modeling type theoretic universes in sheaves, which
some authors have erroneously located in the fact that the amalgamation of a
family of types defined at the leaves of a bar is unique only up to
isomorphism~\citep{xu-escardo:2016:universes}, or in the apparent impredicativity of
certain classic constructions of sheafification~\citep{xu:2015}.

In fact, universes of sheaves are known to be
unproblematic~\citep{streicher:2005}, and sheafification is entirely compatible
with predicative foundations~\citep{awodey-gambino-lumsdaine-warren:2009}; the
real location of the difficulty pointed out by Xu and Escard\'o is that known
constructions of sheafification do not preserve the choice of codes for type
connectives on the nose, so the strict reduction rules that identify
$\mathsf{El}(\hat{A}\hat{\to}\hat{B})$ with
$\mathsf{El}(\hat{A})\to\mathsf{El}(\hat{B})$ must be weakened to canonical
isomorphisms.

These difficulties motivated \citet{coquand-mannaa-ruch:2017} to consider
models of type theory in \emph{stacks}, which express descent from a bar in a
weak enough fashion that universes can be modeled without sheafification.  The
stack model of type theory contains a universe that classifies \emph{discrete}
types (types with no higher-dimensional structure), and has been used to give a
purely denotational account of the independence of Martin-L\"of's type theory
from Markov's principle.

\section*{Acknowledgements}

Thanks to Carlo Angiuli, Mark van Atten, Mark Bickford, Bob Constable, Thierry
Coquand, Mart\'in Escard\'o, Bob Harper, and Per Martin-L\"of for enlightening
conversations about the bar theorem; I am especially thankful to Pierre-Marie
P\'edrot for pointing out an error in a previous draft of this paper, and to
Vincent Rahli for suggesting a significant simplification to the proof, and to
one of the anonymous referees for a suggestion to make the dialogue
interpretation compositional. I am grateful to Bob Constable for scanning and
sending to me a copy of Per Martin-L\"of's doctoral dissertation. In my
formalization, I benefited from Darin Morrison's alternative Agda prelude
library. Finally, I thank the referees for their helpful and constructive
criticisms, and I thank Tristan Nguyen at AFOSR for support.

This work was supported in part by AFOSR under grants MURI FA9550-15-1-0053 and
FA9550-19-1-0216.  Any opinions, findings and conclusions or recommendations
expressed in this material are those of the authors and do not necessarily
reflect the views of the AFOSR.

\section*{Conflicts of Interest}

None.

\nocite{brouwer:1981}

\nocite{fourman:1982,fourman:1984,fourman:2013}
\nocite{vanderhoeven-moerdijk:1984}
\nocite{troelstra-vandalen:1988}
\nocite{brouwer:1981}

\nocite{maclane-moerdijk:1992}
\nocite{sambin:2012}
\nocite{sundholm-van-atten:2008}

\nocite{longley:1999}
\nocite{van-dalen:2013}
\nocite{gambino-schuster:2007}
\nocite{coquand-sambin-smith-valentini:2003}
\nocite{petrakis:2010}
\nocite{capretta-uustalu:2016}
\nocite{constable-bickford:2014}

\nocite{longley-normann:2015}

\nocite{van-heijenoort:2002}

\bibliographystyle{jfplike}
\bibliography{references/refs-bibtex,temp-refs}

\end{document}